\documentclass{amsart}

\let \frak = \mathfrak
\usepackage[mathscr]{euscript}
\let \eusm = \mathscr
\let \bold = \mathbf
\usepackage{amssymb}

\newtheorem{theorem}{Theorem}[section]
\newtheorem{proposition}[theorem]{Proposition}
\newtheorem{lemma}[theorem]{Lemma}
\newtheorem{corollary}[theorem]{Corollary}
\theoremstyle{definition}
\newtheorem*{definition}{Definition}
\theoremstyle{remark}
\newtheorem*{remark}{Remark}
\newtheorem*{remarks}{Remarks}
\newcounter{bean}
\newenvironment{mrk}{\begin{list}{\thebean.}{\usecounter{bean}
\setlength{\leftmargin}{0in} \setlength{\itemindent}{.15in}}}{\end{list}}

\title{The Geometry of Fixed Point Varieties on Affine Flag Manifolds}
\author{Daniel S. Sage}
\address{Department of Mathematics \\
University of Utah \\
Salt Lake City, Utah 84112}
\curraddr{School of Mathematics \\
Institute for Advanced Study \\
Princeton, NJ 08540}
\email{sage@math.ias.edu}

\dedicatory{\large{{\rm To appear in } Transactions of the AMS.}}


\renewcommand{\L}{\Lambda}
\renewcommand{\l}{\lambda}
\renewcommand{\c}{\bold c}
\renewcommand{\b}{\bold b}
\renewcommand{\d}{\delta}
\renewcommand{\a}{\alpha}
\newcommand{\g}{\gamma}
\newcommand{\q}{\bold q}
\newcommand{\Z}{\bold Z}
\newcommand{\N}{\bold N}
\renewcommand{\r}{\bold r}

\newcommand{\Th}{\Theta}
\newcommand{\z}{\zeta}
\newcommand{\s}{\sigma}
\renewcommand{\char}{\text{char}}
\newcommand{\Rs}{R^{Sp}_{sm}}
\newcommand{\Go}{G^{Sp}_{s0}}
\newcommand{\Gs}{G^{Sp}_{sm}}
\newcommand{\RsI}{R^{Sp}_{sI}}
\newcommand{\GsI}{G^{Sp}_{sI}}
\newcommand{\R}{R_{sm}}
\newcommand{\Cs}{C^{Sp}_{sm}}
\newcommand{\C}{C_{sm}}
\newcommand{\CsI}{C^{Sp}_{sI}}
\newcommand{\Ds}{\eusm D^{Sp}}
\newcommand{\Ls}{\eusm L^{Sp}}

\newcommand{\Xs}{X^{Sp}}
\newcommand{\Ys}{Y^{Sp}}
\newcommand{\e}{\epsilon}
\newcommand{\p}{\bold p}
\renewcommand{\d}{\delta}
\renewcommand{\sp}{\frak{sp}_{2n}}
\renewcommand{\Sp}{Sp_{2n}}
\newcommand{\<}{\langle}
\renewcommand{\>}{\rangle}
\newcommand{\m}{\bold m}

\begin{document}
\begin{abstract}

	 Let $G$ be a semisimple, simply connected, algebraic group over an
algebraically closed field $k$ with Lie algebra $\frak g$.  We study the
spaces of parahoric subalgebras of a given type containing a fixed
nil-elliptic element of $\frak g\otimes k((\pi))$, i.e. fixed point
varieties on affine flag manifolds.  We define a natural class of
$k^*$-actions on affine flag manifolds, generalizing actions introduced by
Lusztig and Smelt.  We formulate a condition on a pair $(N,f)$ consisting
of $N\in\frak{g}\otimes k((\pi))$ and a $k^*$-action $f$ of the specified
type which guarantees that $f$ induces an action on the variety of
parahoric subalgebras containing $N$.

	For the special linear and symplectic groups, we characterize all
regular semisimple and nil-elliptic conjugacy classes containing a
representative whose fixed point variety admits such an action.  We then
use these actions to find simple formulas for the Euler characteristics of
those varieties for which the $k^*$-fixed points are finite.  We also
obtain a combinatorial description of the Euler characteristics of the
spaces of parabolic subalgebras containing a given element of certain
nilpotent conjugacy classes of $\frak g$.

\end{abstract}
\maketitle

\section{Introduction}
	Let $G$ be a semisimple, simply connected, algebraic group over an
algebraically closed field $k$, and let $\frak g$ be its Lie algebra.  The
set of Borel subalgebras of $\frak g$ forms a smooth, projective variety
which we denote by $\eusm{B}$.  For any nilpotent element $N$ in $\frak g$,
the set of Borel subalgebras containing $N$ is a closed subvariety
$\eusm{B}_N$ of $\eusm{B}$.  These varieties have been studied extensively
and have interesting applications to representation theory.  For example,
Springer and several others, including Ginzburg, Kazhdan, Lusztig, Slodowy,
Joseph, Borho, and MacPherson, have defined representations of the Weyl
group $W$ in the homology (or cohomology) of the $\eusm{B}_N$'s.  Moreover,
all irreducible representations of $W$ may be obtained from the Springer
representations in the top homology of the $\eusm{B}_N$'s as $N$ runs over
the nilpotent conjugacy classes of $\frak g$.  In fact, when $G=SL_n(k)$,
these representations are precisely the irreducible representations of $W$.

	Kazhdan and Lusztig have introduced the study of the affine
analogue of these fixed point varieties on the flag manifold, which they
anticipate will have applications to the character theory of semisimple
$p$-adic groups.  They have also used affine techniques to extend to all
nilpotent orbits of $\frak g$ a map defined by Carter and Elkington from a
certain class of nilpotent orbits into the conjugacy classes of $W$.  Let
$F=k((\pi))$ be the field of formal power series in one variable with
$A=k[[\pi]]$ its ring of integers.  We then extend scalars to form $G(F)$,
$G(A)$, $\frak g_F=\frak g\otimes {}_{k} F$, and $\frak g_A=\frak g\otimes
{}_{k} A$.  Fix a Borel subalgebra $\frak b \in \eusm{B}$.  An Iwahori
subalgebra of $\frak g_F$ is any $G(F)$-conjugate of the pullback
$\hat{\frak b}$ of $\frak b$ under the natural projection $\frak g_A \to
\frak g$; we denote the space of Iwahori subalgebras by $\hat{\eusm{B}}$.
If $B \subset G$ corresponds to $\frak b$ with preimage $\hat{B} \subset
G(A)$ under the map $G(A) \to G$, then $\hat{\eusm B}$ may be identified
with the quotient space $G(F)/\hat{B}$.  The set $\hat{\eusm{B}}$ has the
structure of an infinite-dimensional projective variety over $k$, i.e. it
is an increasing union of finite-dimensional projective varieties.
(Kazhdan and Lusztig work in the case $k=\mathbf C$, and their proofs of
the results quoted also apply generally in characteristic zero.  However,
for this paper, it will suffice to know that for classical groups in
positive characteristic, these spaces have ind-variety structures and
moreover that certain subspaces are algebraic varieties.  We will show
these results directly.)

	We can also generalize the classical partial flag varieties,
consisting of the subalgebras of $\frak g$ conjugate to a fixed parabolic
subalgebra, to the affine context.  A parahoric subalgebra of $\frak g_F$
is a subalgebra containing an Iwahori subalgebra.  Let $S$ be a set of
simple reflections for the affine Weyl group $\hat{W}$, and for each
$I\subseteq S$, let $\hat{W}_I$ be the subgroup of $\hat{W}$ generated by
$I$.  The parahoric subalgebras containing $\hat{\frak b}$ are precisely
the algebras $\hat{\frak p}_I=\sum_{w\in\hat{W}_I}w\cdot\hat{\frak b}$;
their stabilizers under the adjoint action are the parahoric subgroups
$\hat{P}_I=\hat{B}W_I\hat{B}$.  The map $I\mapsto \hat{\frak p}_I$ gives a
one-to-one correspondence between subsets of $S$ and conjugacy classes of
parahoric subalgebras.  The set $\hat{\eusm{P}}_I$ of $G(F)$-conjugates of
$\hat{\frak p}_I$ is an infinite-dimensional projective variety which may
be identified with the quotient space $G(F)/\hat{P}_I$.  These are the
affine flag manifolds of $G(F)$.  The variety $\hat{\eusm{P}}_S$ is just a
point, so we will only consider the varieties of proper parahoric
subalgebras.

	For each (nontrivial) affine flag manifold $\hat{\eusm{P}}$ and any
element $N$ in the Lie algebra $\frak g_F$, we define $\hat{\eusm{P}}_N$ to
be the space of parahoric subalgebras of the specified type containing $N$.
This is a closed subvariety of $\hat{\eusm{P}}$.  We shall restrict our
attention to the topologically nilpotent (or nil) elements $N$, those elements
such that $(\text{ad}\, N)^m \to 0$ in the power series topology on End
$\frak g_F$ as $m \to \infty$.  More concretely, these are the elements of
$\frak g_F$ which may be conjugated into some $N^\prime \in \frak g_A$
whose image under projection to $\frak g$ is nilpotent.  Kazhdan and
Lusztig have shown in \cite{KL} that if $\frak g$ is simple, then for $N$ nil, the closed subvariety
$\hat{\eusm{P}}_N$ is finite-dimensional (although possibly having
infinitely many irreducible components) precisely for $N$ regular
semisimple.  If we further assume that $N$ is elliptic, i.e. the connected
centralizer of $N$ in $G(F)$ is an anisotropic maximal torus, then
$\hat{\eusm{P}}_N$ is a projective variety in the usual sense.

	The affine fixed point varieties $\hat{\eusm{P}}_N$ are much more
complicated than their classical counterparts, and not much is known about
their topology.  Unlike the $P_N$'s, they need not be rational nor do they
necessarily have cell decompositions.  For example, even in the low
dimensional case $\frak {sp}_6(F)$, Bernstein and Kazhdan have found $N$
such that an irreducible component of $\hat{\eusm{B}}_N$ admits a dominant
morphism onto an elliptic curve \cite{KL}.  In fact, even the dimensions
of these varieties were not known in general until recently, when
Bezrukavnikov proved a conjectural formula of  Kazhdan and Lusztig.

	Lusztig and Smelt have defined $k^*$-actions on certain of these
fixed point varieties for the special linear groups \cite{LS}.  They used these
actions to obtain interesting formulas for the Euler characteristics of
these subvarieties for the varieties of maximal parahoric and Iwahori
subalgebras.  We have been studying
those $\hat{\eusm{P}}_N$ which admit a class of algebraic $k^*$-actions
which generalize the Lusztig-Smelt actions.  In
section 2, we define a natural action of the group
$Q'_G=G_{\text{ad}}(k)\times k^*$ on the affine flag manifolds of $G(F)$.
An algebraic homomorphism $k^*\to Q'_G$ defines a $k^*$-action on
$\hat{\eusm P}$ which we call a $Q'_G$-action.  The group $Q'_G$ also acts
on $\frak g_F$ by $(g,\lambda)\cdot N=\text{Ad}(g)N^\lambda$ for $N\in\frak
g_F$, where $N^\lambda$ is the image of $N$ under the $k^*$-action
on $\frak{g}_F$ induced by mapping the uniformizing parameter $\pi$ of $F$ to $\lambda\pi$.  We say that a $k^*$-action $f$ almost commutes with an element $N\in
\frak g_F$ if for some nonzero integer $s$, $N$ is an eigenvector of
$f(\lambda)$ with eigenvalue $\lambda^s$ for all $\lambda\in k^*$.  We show
that if $(N,f)$ is an almost commuting pair, then $f$ induces an action on
$\hat{\eusm P}_N$.

	In order to study these actions, we need a more concrete
description of our infinite-dimensional varieties, analogous to the
identification of the varieties of parabolic subalgebras of $\frak g$ with
certain varieties of partial flags.  Let $\frak g$ be a classical Lie
algebra, and take the standard representation $\frak g \hookrightarrow
\frak {gl}(V)$.  Set $n$ equal to the dimension of $V$.  The varieties of
maximal parahoric subalgebras simply consist of certain classes of
$A$-lattices in $V_F$.  The full affine flag manifold $\hat{\eusm{B}}$ is a
space of complete chains of lattices ${L_0 \supset L_1 \supset \dots
\supset L_n}$, where $L_i$ is of $k$-codimension $1$ in $L_{i-1}$ for each
$i$, $L_n = \pi L_0$, and the lattices $L_i$ have certain properties
depending on $\frak g$.  The $\hat{\eusm{P}}$ for other parahorics are sets
of partial lattice chains.  Since the Lie algebra action of $\frak g_F$ on
$\hat{\eusm{P}}$ corresponds to the natural action of $\frak g_F$ on
lattice chains, the variety $\hat{\eusm{P}}_N$ may be identified with those
partial lattice chains fixed by $N$.  We can now define a $k^*$-action on
$\hat{\eusm{P}}_N$ by defining a $k^*$-action on $V_F$ which preserves
lattices (and hence lattice chains) of the specified types and then
descends to the fixed point subvarieties.

	In sections 3 and 4, we identify the affine flag manifolds for
$\frak{sl}(V_F)$ and $\frak{sp}(V_F)$ with varieties of lattice chains. We
then apply the results of section 2 to obtain complete and concrete
descriptions, in terms of the characteristic polynomial, of the
nil-elliptic conjugacy classes containing a representative $N$ whose fixed
point varieties admit $k^*$-actions induced by a $Q'_G$-action. (For the
symplectic algebras, we assume that the characteristic of $k$ is not $2$.)
Moreover, for each such conjugacy class, we find class representatives in a
uniform way for which the lattices fixed by both $N$ and the $k^*$-action
can be determined explicitly.

	These actions are useful and interesting because they enable one to
determine the Euler characteristics of the $\hat{\eusm{P}}_N$'s.  By a
theorem of Bialynicki-Birula \cite{B-B1}, if the multiplicative group $k^*$
acts on an algebraic variety $Y$, then the Euler characteristics of $Y$ and
$Y^{k^*}$ are the same, where $Y^{k^*}$ denotes the fixed point set of the
action.  Therefore, we need only understand the much simpler variety of
$k^*$-fixed partial lattice chains to calculate the Euler characteristics
of the $\hat{\eusm{P}}_N$'s.  In the most favorable cases, including every
nil-elliptic conjugacy class of $\frak{sl}(V_F)$ obtained in the above
classification, these $k^*$-actions only fix diagonal lattices relative to
a fixed $k$-basis $e_1,\dots,e_n$ of $V$, i.e. lattices of the form
$\eusm{L}_{\bold r} = A\pi ^{r_1}e_1+\dots +A\pi ^{r_n}e_n$ for $\bold r
\in \bold Z ^n$.  This implies that the fixed point set
$(\hat{\eusm{P}}_N)^{k^*}$ is discrete and in fact finite, since $N$ is
nil-elliptic.  We have thus reduced the computation of the Euler
characteristic of $\hat{\eusm{P}}_N$ to a purely combinatorial problem of
counting lattice chains.  For those pairs $(N,f)$ such that $f$ fixes only
diagonal lattices, we have discovered simple formulas for the Euler
characteristics of $\hat{\eusm{P}}_N$ for each class of parahoric
subalgebras as products of binomial coefficients.  These formulas can be
interpreted as the number of simplices of type corresponding to the
parahoric class within a certain  simplex in an apartment of the affine
building of $G$.  It is possible to apply
similar techniques to the orthogonal algebras; this will be the subject of
a future paper.

If an algebraic $k^*$-action on a
smooth projective variety $Y$ has isolated fixed
points, the cells of the Bialynicki-Birula plus decomposition of $Y$ are
isomorphic to affine space \cite{B-B2}.  The cells are indexed by the fixed
points of the action, and the cell corresponding to a fixed point $x$ is
$\{y\in Y \mid \lim_{\lambda\to 0}\lambda\cdot y =x\}$, where the limit is
understood in the sense of \cite{B-B3}.  It is natural to conjecture that
when the fixed point set $(\hat{\eusm{P}}_N)^{k^*}$ is finite, the variety
$\hat{\eusm{P}}_N$ has such a
cell decomposition.  For $\frak{sl}(V)$, an argument used by Lusztig and
Smelt to obtain this result for Iwahori and maximal parahoric subalgebras applies
to show that this is indeed true.  In particular, the \'{e}tale cohomology
of $\hat{\eusm{P}}_N$ vanishes in odd degrees in these cases.

	These $k^*$-actions also have applications to the topology of the
classical fixed point varieties.  This is in very much the same spirit as
Kazhdan and Lusztig's use of nil-elliptic classes in $\frak g_F$ to study
nilpotent classes in $\frak g$.  Let $X$ denote the variety of conjugates
of the maximal parahoric subalgebra $\frak g_A$.  The natural
$\eusm{B}(k)$-fibration $\hat{\eusm{B}} \overset {\rho}{\to} X$ is given in
terms of lattice chains by $(L_0,\dots,L_n) \overset{\rho}{\mapsto} L_0$.
The fibers of this map are the complete flags (of the type corresponding to
$\frak g$) in the $k$-vector space $L_0/\pi L_0$, and the $N$-fixed lattice
chains in the fiber may be interpreted as $\eusm{B}_{\bar{N}}$, where
$\bar{N}$ is the induced nilpotent transformation on $L_0/\pi L_0$.  If
$(\hat{\eusm{B}}_N)^{k^*}$ is finite and $L_0\in (X_N)^{k^*}$, then the
Euler characteristic of $\eusm{B}_{\bar{N}}$ is the number of lattice
chains in $(\hat{\eusm{B}}_N)^{k^*}$ lying above $L_0$.  Similar
considerations apply to the partial flag varieties.  In section 3, we show
that every nilpotent conjugacy class of $\frak{sl}(V)$ may be realized as
the class of such an $\bar{N}$.  We thus obtain a combinatorial
description of the Euler characteristics of the classical fixed point
varieties for $\frak{sl}(V)$ without using their explicit cell
decomposition or arguments involving induction on unipotent classes.  We
use a similar technique in section 4 to find combinatorial interpretations
of the Euler characteristics of the classical fixed point varieties for
$\frak{sp}(V)$ for certain nilpotent classes.

	This paper is based on part of my doctoral thesis at the University
of Chicago \cite{S}.  It is a great pleasure to thank my advisor, Robert
Kottwitz, for many extremely helpful comments and suggestions.  I am also
happy to acknowledge the support of a National Science Foundation Graduate
Fellowship and a University of Chicago McCormick Fellowship.

\section{$\lowercase{k}^*$-actions on fixed point varieties of affine flag
manifolds}

	Let $k$ be an algebraically closed field, and let $G=G(k)$ be a
semisimple, simply connected, algebraic group over $k$ with Lie algebra
$\frak g$.  Let $F=k((\pi))$ be the field of formal power series in one
variable with ring of integers $A=k[[\pi]]$.  The local field $F$ is
furnished with its usual valuation so that $v(\pi)=1$.  We now extend scalars
of the group and Lie algebra to the power series field and ring.

	Let $p:\frak {g}_A\to \frak {g}$ be the natural projection.  An
Iwahori (or minimal parahoric) subalgebra of $\frak g_F$ is a
$G(F)$-conjugate of the pullback via $p$ of any Borel subalgebra of
$\frak{g}$.  The parahoric subalgebras of $\frak g_F$ are those
subalgebras which contain an Iwahori subalgebra.  Iwahori and parahoric
subgroups of $G(F)$ may be defined in the same manner.  Alternatively, a
parahoric subgroup is the stabilizer of a parahoric subalgebra under the
adjoint action.  If $\hat{\frak p}$ is a parahoric subalgebra with
stabilizer $\hat{P}$, then the set of parahoric subalgebras conjugate to
$\hat{\frak p}$ may be identified with the quotient space $G(F)/\hat{P}$.
This space has the natural structure of an infinite-dimensional projective
algebraic variety over $k$, i.e. it is an increasing union of
finite-dimensional projective varieties.  (The results on the variety structure
of these spaces are proved by Kazhdan and Lusztig for $k=\mathbf C$ (and
their proofs work generally in characteristic zero).  For classical groups,
we may give these spaces the structure of ind-varieties, using their
concrete identifications with lattices and lattice chains as explained in
section three.  For our computations of the Euler characteristic of
certain fixed point subvarieties, it will suffice to know that the spaces
in question are indeed algebraic varieties in positive characteristic.
Since we will show this directly, in the present paper we will not seek to
extend the general results of Kazhdan and Lusztig to characteristic $p$.)

	Fix a Borel subalgebra $\frak b$ of $\frak g$ with $B$ the
corresponding Borel subgroup, and let $T$ be a maximal torus of $G$
contained in $B$.  The affine Weyl group is a Coxeter group defined by
$\hat{W}=N(T(F))/T(A)$, where $N(T(F))$ is the normalizer of $T(F)$ in
$G(F)$.  Let $S$ be a set of simple reflections for $\hat{W}$.  Then the
conjugacy classes of parahoric subalgebras correspond bijectively to the
subsets of the simple reflections $S$.  Given $I\subseteq S$, let
$\hat{W}_I$ be the subgroup of $\hat{W}$ generated by $I$.  We obtain a
standard Iwahori subalgebra $\hat{\frak b}$ (with respect to the choice of
$\frak b$) by pulling back $\frak b$ to $\frak g_A$.  There is precisely
one parahoric subalgebra from each conjugacy class containing $\hat{\frak
b}$; these are just the subalgebras $\hat{\frak
p}_I=\sum_{w\in\hat{W}_I}w\cdot\hat{\frak b}$, and with $\hat{B}$ the
stabilizer of $\hat{\frak b}$, their stabilizers are
$\hat{P}_I=\hat{B}\hat{W}_I\hat{B}$.  The corresponding varieties of
conjugate subalgebras will be denoted $\hat{\eusm P}_I$.  If $I\subseteq
J$, then $\hat{P}_I\subseteq\hat{P}_J$, inducing a natural map $\hat{\eusm
P}_I\to\hat{\eusm P}_J$; this is a fiber bundle with a classical (partial)
flag variety as fiber.  Note that $\hat{\eusm P}_{S}$ is just a point.
From now on, we shall restrict attention to the nontrivial affine flag
manifolds or equivalently, to the proper parahoric subalgebras.

	We are interested in defining actions of $k^*$ on the subvarieties
of the affine flag manifolds containing a fixed element $N$ of $\frak g_F$.
As a first step, we will define a natural class of $k^*$-actions on the
affine flag manifolds themselves.  The multiplicative group $k^*$ acts by
$k$-algebra automorphisms on $F$ via the formula
$\lambda\cdot\pi^m=\lambda^m\pi^m$ for all $\lambda\in k^*$ and $m\in\bold
Z$.  This induces an action of $k^*$ on $G(F)$.  Since $k^*$ preserves
the standard Iwahori subgroup $\hat{B}$ and acts trivially on
the affine Weyl group, the action descends to the varieties of
parahoric subalgebras.  Also, $G(F)$ acts on itself by left
multiplication.  It is now easy to see that the semidirect product
$G(F)\rtimes k^*$ acts on $G(F)$ by $(g,\lambda)\cdot h=gh^\lambda$ where
$g,h\in G(F)$ and $\lambda\in k^*$ and induces an action on the
affine flag manifolds.  Since the center $Z(F)$ of $G(F)$ is contained in
every parahoric subgroup, we finally obtain an action of
$Q_G=G_{\text{ad}}(F)\rtimes k^*$ on each $\hat{\eusm P_I}$, where
$G_{\text{ad}}$ is the adjoint form of $G$.  We have
$G_{\text{ad}}(F)^{k^*}=G_{\text{ad}}(k)$, so the group $Q_G$ contains the
direct product $Q'_G=G_{\text{ad}}(k)\times k^*$.  We can now define the
class of $k^*$-actions which we will study.

\begin{definition} A $k^*$-action $f$ on an affine flag
manifold for the group $G$ is called a $Q'_G$-{\it action} if it is
obtained from an algebraic homomorphism $k^*\overset{f}{\to} Q'_G$.
\end{definition}

\begin{remark} We can similarly define $Q_G$-actions, but in
this paper, we will restrict attention to $Q'_G$-actions.
\end{remark} 
	
	The group $Q_G$ also acts on $\frak g_F$ by $(g,\lambda)\cdot
N=\text{Ad}(g)N^\lambda$ for $N\in\frak g_F$.

\begin{definition} Given $N\in \frak g_F$, we say that
a $k^*$-action $f$ {\it almost commutes} with $N$ if for some nonzero
integer $s$, $N$ is an eigenvector of $f(\lambda)$ with eigenvalue $\lambda^s$,
i.e. $f(\lambda)\cdot N=\lambda^s N$ for all $\lambda\in k^*$.  The integer $s$
is called the exponent of the pair $(N,f)$.
\end{definition}

This property guarantees that a $Q_G$-action $f$ on a conjugacy class of
parahoric subalgebras $\hat{\eusm P}$ descends to $\hat{\eusm P}_N$, the
subvariety of parahorics containing $N$.

\begin{lemma} Let $\hat{\eusm P}$ be an affine flag manifold, and suppose
that $N$ almost commutes with the $Q_G$-action $f$.  Then $f$ induces an
action on $\hat{\eusm P}_N$.
\end{lemma}

\begin{proof} Take $\hat{\frak p}\in\hat{\eusm P}_N$, and set
$f(\l)=(g(\l),\l^p)$ for some $p\in\Z$.  (The map of $k^*$ into the second
factor of $Q_G$ is a homomorphism from $k^*$ to itself.) Then
$N=\l^{-s}f(\l)\cdot N=\text{Ad}\, g(\l)N^{\l^p}\in\text{Ad}\, g(\l)\hat{\frak
p}^{\l^p}=f(\l)\cdot\hat{\frak p}$ as desired.
\end{proof}

	Elements of $\frak g_F$ which almost commute with a $Q_G$-action
must satisfy the following condition:

\begin{theorem} Let $\rho:G\to GL(V)$ be a representation of $G$,
inducing the Lie algebra representation $d\rho:\frak g_F\to\frak
{gl}(V_F)$.  Suppose that $N$ almost commutes with the $Q_G$-action
$f(\l)=(g(\l),\l^p)$ with exponent $s$.  Then the characteristic polynomial
$\text{char}_{d\rho(N)}(\mu)$ is a homogeneous polynomial of degree
$n=\dim\ V$ in $\mu$ and $\pi^q$ for some rational number $q\in
\bigcup_{i=1}^n\frac1i\Z$; if $d\rho(N)$ is not nilpotent, then $p\ne 0$ and
$q=s/p$.
\end{theorem}

\begin{proof} Set $\text{char}_{d\rho(N)}(\mu)=\sum_{i=0}^n c_i \mu^{n-i}$ for some $c_i\in
F$.  Since $d\rho(\text{Ad}\, g(\l)N^{\l^p})$ and $d\rho(N^{\l^p})$ are
$GL(V_F)$-conjugate, we see that $\sum_{i=0}^n c_i^{\l^p}
\mu^{n-i}=(\text{char}_{d\rho(N)}(\mu))^{\l^p}=\text{char}_{d\rho(N^{\l^p})}(\mu)=\text{char}_{d\rho(\l^sN)}(\mu)=\sum_{i=0}^n
\l^{si}c_i \mu^{n-i}$.  The claim is trivial if $d\rho(N)$ is nilpotent, so
we may assume that $c_j\ne 0$ for some $j\ge 1$.  The equations
$c_i^{\l^p}=\l^{si}c_i$ imply that $pv(c_i)=si$ for each $i$ with $c_i\ne
0$; furthermore, $c_i$ is homogeneous, i.e. a Laurent series in $\pi$ with
at most one nonvanishing term.  This shows that $\text{char}_{d\rho(N)}(\mu)$ is homogeneous of degree
$n$ in $\mu$ and $\pi^q$, where $q=s/p$, and since $sj/p$ is an integer,
$q\in\frac1j\Z$.
\end{proof}

\begin{definition} If $N$ is not nilpotent, we call the
ratio $s/p$ the {\it loop index} of the pair $(N,f)$.
\end{definition}

\begin{corollary} Let $\overline F$ be an algebraic closure of $F$.  Then
all nonzero eigenvalues of $d\rho(N)$ are homogeneous elements of
$\overline{F}$ with valuation equal to the loop index.
\end{corollary}

\begin{proof} By the theorem, $\text{char}_{d\rho(N)}(\mu)=\sum_{i=0}^n c_i
(\pi^q)^i\mu^{n-i}$ for some $c_i\in k$.  The roots of $\text{char}_{d\rho(N)}(\mu)$
in the algebraic closure $\overline{F}$ are just $x\pi^q$ where $x\in k$
is a root of the polynomial $\sum_{i=0}^n c_i\mu^{n-i}\in k[\mu]$.
\end{proof}

	Recall that $N\in\frak g_F$ is called {\it nil} (or
topologically nilpotent) if $\lim_{m\to\infty}(\text{ad }N)^m=0$.  We are
primarily interested in defining $k^*$-actions on the fixed point variety
$\hat{\eusm P}_N$ when $N$ is regular semisimple and nil, since in this case,
$\hat{\eusm P}_N$ is  finite-dimensional \cite{KL}.  For $N$ of this type in the
case of the special linear and symplectic Lie algebras, we will
prove a converse statement to theorem 2.2, namely, that if the
charactersitic polynomial of a nil, regular semisimple conjugacy class
has the form given in theorem 2.2, then we can find an element of the class
together with a $Q_G$-action almost commuting with it. 

	Now assume that $G$ is a classical group, and let $\rho:G\to GL(V)$
be the standard representation of $G$, inducing the Lie algebra
monomorphism $\frak g\hookrightarrow \frak{gl}(V)$.  In order to define
explicit $Q'_G$-actions when $G$ is a classical group, it will be
convenient to consider homomorphisms $k^*\to
Q'\overset{\text{def}}{=}GL(V)\times k^*$.  Such a homomorphism gives rise
to a $Q'_{SL(V)}$-action via the projection $GL(V)\to PGL(V)=PSL(V)$, and if
$G$ is not the special linear group, we can obtain a $Q'_G$-action in a
similar way after placing suitable restrictions on the homomorphism.

	More generally, if $V$ is an $n$-dimensional vector space over $k$,
$Q'$ acts on $V_F$ by $k$-linear automorphisms.  We call a $k^*$-action $f$
on $V_F$ a $Q'$-{\it action} if it comes from an algebraic homomorphism
$k^*\overset{f}{\to} Q'$.  For $n\ge 2$, this is equivalent to $f$ mapping
$A$-submodules of $V_F$ to $A$-submodules and the existence of a $k$-basis
$\underline e=(e_1,\dots,e_n)$ of $V$ such that every element of the
topological spanning set $\underline{\hat{e}}=\{\pi^m e_i|m\in\bold Z,
i\in[1,n]\}$ of $V_F$ is an eigenvector of $f(\lambda)$ for each
$\lambda\in k^*$.  In fact, the eigenvalue of $\pi^m e_i$ must have the form
$\l^{\nu(\cdot,i)}$, where $\nu(\cdot,i):\bold Z\to\bold Z$ is an affine
function with slope independent of $i$, i.e. $\nu(m,i) = pm+c_i$ for some
integers $p$ and $c_i$ \cite[theorem 2.5]{S}.  (A specific example of a
$k^*$-action of this form was first introduced by Lusztig and Smelt
\cite{LS}.)

	We say that an endomorphism $N\in\frak{gl}(V_F)$ almost commutes
with $f$ if for some nonzero integer $s$, $f(\lambda) N=\lambda^s N
f(\lambda)$ for all $\lambda\in k^*$.  The group $Q'$ acts on
$\frak{gl}(V_F)$, and noting that $f(\l)\cdot N=f(\lambda)
Nf(\lambda)^{-1}$, we see that this is compatible with our previous
definition.  Fix a basis $\underline e$ such that the elements of
$\underline{\hat{e}}$ are eigenvectors of the $f(\l)$'s.  Let
$(a_{ij})\in\frak{gl}_n(F)$ be the matrix for $N$ determined by $\underline
e$, and write $a_{ij}=\sum_{l\gg -\infty} a_{ij}^l \pi^l$ where
$a_{ij}^l\in k$.  Equating the $\pi^{m+l}e_j$ terms in $f(\lambda) N \pi^m
e_i=\lambda^s N f(\lambda) \pi^m e_i$ for all $m$ and $i$ gives the
transformation laws $ \nu(m,j)+s = \nu(m+l,i)$ for all $m\in\bold Z$ and
all $i$, $j$, and $l$ such that $a_{ij}^l\ne 0$.

	The characteristic polynomial determines the $GL(V_F)$-conjugacy
classes of semi\-simple endomorphisms.  We can now characterize the regular
semisimple conjugacy classes obtained from almost commuting pairs $(N,f)$
with $f$ a $Q'$-action.

\begin{theorem} Let $\eusm C$ be a regular semisimple conjugacy class in
$\frak{gl}(V_F)$.  Then there is a representative $N\in \eusm C$ almost
commuting with a $Q'$-action if and only if the characteristic polynomial
$\text{char}_{\eusm C}(\mu)$ is a homogeneous polynomial of degree $n$ in
$\mu$ and $\pi^{q}$ for some nonzero rational number $q\in
\frac1{n-1}\Z\cup \frac1n \Z$.
\end{theorem}

\begin{proof} The forward implication follows just as in theorem 2.2.  We
need only note that since at most one of the eigenvalues is $0$, either the
$\pi^{qn}\mu^0$ or the $\pi^{q(n-1)}\mu$ term of the characteristic
polynomial has nonzero coefficient, implying that $q\in \frac1{n-1}\Z\cup
\frac1n \Z$.

	Conversely, assume that for some nonzero $q\in \frac1{n-1}\Z\cup \frac1n
\Z$, we have $\text{char}_{\eusm C}(\mu)=\sum_{i=0}^n
c_i(\pi^q)^i\mu^{n-i}$ with each $c_i\in k$ and $c_0=1$.  Define a matrix $N=(a_{ij})$
by
\[
a_{ij}= \begin{cases} 1& \text{if $j=i+1$}\\ 
(-1)c_i\pi^{iq}& \text{if $j=1$}\\ 
0& \text{otherwise.} \end{cases} \] 
Note that if $c_i\ne 0$, then $qi\in \bold Z$.  Thus, $a_{ij}\in F$ for
each $i$ and $j$.

	If the determinant of $N$ is nonzero, set $\nu(m,i)=nm-inq$;
otherwise, set $\nu(m,i)=(n-1)m-i(n-1)q$.  It is easy to check that if we
define $f(\l)\pi^m e_i=\l^{\nu(m,i)}\pi^m e_i$, where $(e_1,\dots,e_n)$ is
the standard basis of $k^n$, then $N$ and $f$ almost commute with exponent
$nq$ (resp. $(n-1)q$).  (The existence of a $Q'$-action almost commuting
with $N$  is also guaranteed by theorem
2.20 of \cite{S}.)

	Finally, it follows easily from the definition of the determinant
as a sum over the symmetric group $S_n$ that
$\text{char}_N=\text{char}_{\eusm C}$.  The conjugacy class $\eusm C$ is
regular semisimple, so $\text{char}_{\eusm C}$ has $n$ distinct roots, and
since any endomorphism in $\frak{gl}(V_F)$ with $n$ distinct eigenvalues must
be diagonalizable, we have $N\in \eusm C$.
\end{proof}

	It is known that an endomorphism is nil if and only if there exists
$N'\in\frak{gl}(V_A)$ conjugate to $N$ such that the image of $N'$ under
projection to $\frak{gl}(V)$ is nilpotent \cite{KL}.  It is immediate that
all nonzero eigenvalues of a nil endomorphism must have positive valuation.
Since a regular semisimple endomorphism with characteristic polynomial as
above has eigenvalues of the form $x\pi^q$ for some $x\in k$, we obtain the
following corollary:

\begin{corollary} Let $\eusm C$ be a nil, regular semisimple conjugacy
class in $\frak{gl}(V_F)$.  Then there is a representative $N\in \eusm C$
almost commuting with a $Q'$-action if and only if the characteristic
polynomial $\text{char}_{\eusm C}(\mu)$ is a homogeneous polynomial of
degree $n$ in $\mu$ and $\pi^{q}$ for some positive rational number
$q\in \frac1{n-1}\Z\cup \frac1n
\Z$.  
\end{corollary}

\section{$\lowercase{k}^*$-actions on fixed point varieties of affine flag
manifolds: The special linear algebras}

	We will now apply our results on $k^*$-actions to the study of
fixed point varieties on affine flag manifolds.  In this section, we will
consider the affine flag manifolds obtained for the root systems $A_n$.

	Let $V$ be a $k$-vector space of dimension $n\ge2$.  The special
linear group $SL(V)$ is the group of linear automorphisms of $V$ which
induce the identity on the highest exterior power of $V$, i.e. the volume-preserving linear
automorphisms.  Its Lie algebra is denoted $\frak{sl}(V)$.  Fixing a basis
$e_1,\dots,e_n$ identifies $SL(V)$ with $SL_n(k)$, the semisimple, simply
connected, algebraic group over $k$ of type $A_{n-1}$.

	Let $\frak b\subset \frak{sl}_n(k)$ be the standard Borel subalgebra
determined by the fixed basis $\underline{e}$, namely the upper triangular matrices in $\frak
{sl}_n(k)$.  The standard Iwahori algebra $\hat{\frak b}$ is the pullback of $\frak
b$ to $\frak{sl}_n
(A)$.  We now define algebras $\hat{\frak q}_i\subset\frak{sl}_n(F)$ for $i\in[0,n-1]$ by 
\begin{align*}
 &\quad\qquad\begin{matrix} \negthickspace\negthickspace i&\thickspace\thickspace n-i\end{matrix}\\ 
\hat{\frak q}_i=\frak{sl}_n(F) &\cap\begin{pmatrix}
A&\pi^{-1} A\\ \pi A&A \end{pmatrix} \begin{matrix} i\\n-i.
\end{matrix}
\end{align*}   
The notation means that $\hat{\frak q}_i$ consists of those matrices in
$\frak{sl}_n(F)$ of the above form, i.e. the upper right block is
an $i\times i$ matrix with entries in $A$ and similarly for the other
blocks.  Blocks with zero as a dimension are ignored, so $\hat{\frak
q}_0=\frak{sl}_n(A)$.  These are the maximal (proper)
parahoric subalgebras containing $\hat{\frak b}$.  The stabilizer of
$\hat{\frak q}_i$ is the group $\hat{Q}_i$ of matrices in $SL_n(F)$
with the same block structure.  We can now obtain all parahorics containing
$\hat{\frak b}$ in terms of the maximal ones: $\hat{\frak p}_I=\bigcap_{i\notin
I} \hat{\frak q}_i$ with stabilizer $\hat{P}_I=\bigcap_{i\notin I}
\hat{Q}_i$.  It is possible to choose simple reflections for the
affine Weyl group $\hat W$, indexed by the set $[0,n-1]$, such that this
notation agrees with the notation of section 2 \cite{S}.

	The affine flag manifolds of $SL(V_F)$ may be interpreted as
varieties of chains of $A$-lattices in $V_F$.  Let $\eusm L$ denote the
space of all $A$-lattices in $V_F$.  This set has a natural ind-variety
structure.  To see this, note that for any lattice $L$, there exists a
natural number $m$ such that $\pi^mV_A\subseteq L\subseteq \pi^{-m}V_A$.
The natural map $L\mapsto L/\pi^mV_A$ from the set of lattices satisfying this condition
for a fixed $m$ to the variety of subspaces of $\pi^{-m}V_A/\pi^mV_A$ is an
embedding, identifying the domain as a projective variety.  Taking the
union over $m$ makes $\eusm L$ into an infinite-dimensional variety.
(For classical groups, affine flag manifolds can be interpreted in terms of
lattices, and similar arguments show that they are ind-varieties.)

If $L$ is a lattice with basis
$f_1,\dots,f_n$, then $f_1\wedge\dots\wedge f_n=\alpha e_1\wedge\dots e_n$
for some $\alpha\in F^*$.  The valuation of $\alpha$ is independent of the
choice of bases for $L$ and $V$; we call it the {\it valuation}  of
$L$ and denote it by $v(L)$.  The space of lattices with fixed valuation $m$ is denoted $\eusm
L_m$. The group $GL(V_F)$ acts transitively on $\eusm L$, and it is obvious
that $T\in GL(V_F)$ maps $\eusm L_m$ to $\eusm L_{m+v(\det T)}$.  The
$SL(V_F)$-orbits of $\eusm L$ are just the varieties $\eusm L_m$ for $m\in\bold
Z$.

	The actions of $SL(V_F)$ and $F^*$ on $\eusm L$ commute, so we see
that $\eusm L_m\cong\eusm L_{m'}$ as an $SL(V_F)$-variety when $m\equiv
m'\pmod n$.  (All the $\eusm L_m$'s are isomorphic as varieties, but the
$SL(V_F)$-structures are different.)  As a result, for each $i\in [0,n-1]$,
$SL(V_F)$ acts on the spaces $X_i=\{\dots\supset\pi^{-1}L\supset
L\supset\pi L\supset\dots\mid L\in\eusm L_i\}$.  The set $X_i$ is just the
variety of lattices modulo homothety with valuation congruent to $i$ modulo
$n$.  This alternate point of view shows that the $\eusm L_m$'s are just
special cases of the following more general class of varieties:

\begin{definition} Let $I$ be a nonempty subset of $[0,n-1]$.  We call
\[X_I=\{(L_i)\mid i\in I\negthickspace\negthickspace\pmod n,\, L_i\in\eusm L_i,\, L_i\supset L_j\iff
i<j,\text{ and $L_{i+n}=\pi L_i$}\}\] 
the space of (partial) lattice chains
of type $I$.  
\end{definition}

	The varieties $X_I$ are in fact nothing more than the affine flag
manifolds for $SL(V_F)$.  Note that if $I\subseteq[0,n-1]$ contains $0$,
then for any $L\in\eusm L_0$, the fiber over $L$ of the projection $X_I\to
X_{\{0\}}=\eusm L_0$ is isomorphic as an $SL(L/\pi L)$-variety to $\eusm
F_{I\setminus\{0\}}$, the variety of flags of type $I\setminus\{0\}$ in
$L/\pi L$.  Using this observation, we can show that $SL(V_F)$ acts transitively on
each $X_I$.  The automorphism of $\Z/n\Z$ given by $c\mapsto -c$ defines an
involution $\alpha$ on the set of coset representatives $[0,n-1]$.  The map
$I\mapsto\alpha(I)$ is thus a bijective self-map on the set of subsets of
$[0,n-1]$.  We can find a lattice chain in $X_I$ whose stabilizer is
$\hat{P}_{\alpha(I)^c}$.  This proves the following theorem (for more
details, see \cite[theorem 3.4]{S}):

\begin{theorem} The varieties $X_I$ and $\hat{\eusm P}_{\alpha(I)^c}$ are
isomorphic for each nonempty subset $I$ of $[0,n-1]$.  
\end{theorem}

	In order to apply the results of section 2 on almost commuting
pairs to the present context, we will relate regular semisimple conjugacy
classes in $\frak{sl}(V_F)$ and $\frak{gl}(V_F)$.  Fix an algebraic closure
$\overline F$ of $F$, and let $F_s$ be the separable closure of $F$ in
$\overline F$.  Recall that for the simple root systems, the torsion primes
are $2$ for $B_n$, $D_n$, and $G_2$, $2$ and $3$ for $E_6$, $E_7$, and
$F_4$, and $2$, $3$, and $5$ for $E_8$ \cite[I.4.4]{SS}.

\begin{proposition} Let $G$ be a simple, simply connected, algebraic
group over $k$ with Lie algebra $\frak g$, and assume that the
characteristic of $k$ is
not a torsion prime for $G$.  If $N$ is a regular semisimple element of
$\frak g_F$, then the $G(F)$-conjugacy class of $N$ in $\frak g_F$ is the
intersection of $\frak g_F$ with the $G(\overline F)$-conjugacy class of
$N$ in $\frak g_{\overline F}$.  
\end{proposition}

\begin{proof} Let $N'\in\frak g_F$ be a $G(\overline F)$-conjugate of $N$.
We claim that $N$ and $N'$ are $G(F_s)$-conjugate.  Let $Z=Z(\overline F)$
and $Z'=Z'(\overline F)$ be the centralizers of $N$ and $N'$ respectively
in $G(\overline F)$. Since the endomorphisms are semisimple, the
corresponding Lie algebras $\frak z$ and $\frak z'$ are the centralizers of
$N$ and $N'$ in $\frak g_{\overline F}$ \cite[p. 128]{Bo}.  The group $G$ is simply connected
and the characteristic of $F$ is not a torsion prime, so by a theorem of Springer and
Steinberg, $Z$ and $Z'$ are connected \cite[II.3.19]{SS}.  Since
$N$ and $N'$ are regular semisimple, they are maximal tori defined over
$F$.  The tori $Z$ and $Z'$ split over $F_s$ \cite[III.8.11]{Bo}, hence are
conjugate over $G(F_s)$.  Choose $h\in
G(F_s)$ such that $hZh^{-1}=Z'$, implying that $h\frak zh^{-1}=\frak z'$.
But $G(\overline F)$-conjugacy in the Cartan subalgebra $\frak z'$ is
equivalent to $W$-conjugacy, where $W$ is the Weyl group.  Since
$W=N(Z'(F_s))/Z'(F_s)$, there exists $h'\in N(Z'(F_s))$ such that
$h'hNh^{-1}h'{}^{-1}=N'$ as desired.

	Set $g=h'h$.  For each $\gamma\in\text{Gal}(F_s/F)$, we have
$\g(g)N\g(g^{-1})=N'$, so that $g^{-1}\g(g)$ is in $Z(F_s)$.  The map
$\g\mapsto\g^{-1}\g(g)$ is a cocycle with coefficients in the torus
$Z(F_s)$; we show that it is a coboundary.  The field $F$ is complete with
respect to a discrete valuation, and its residue field is algebraically
closed.  Therefore, the Brauer group of every finite separable extension of $F$
is zero, and this implies that $H^1(\text{Gal}(F_s/F),T(F_s))=0$ for any
torus $T$ defined over $F$ \cite[p. 170]{S1}.  In particular,
$H^1(\text{Gal}(F_s/F),Z(F_s))=0$.  Thus, there exists
$z\in Z(F_s)$ such that $z^{-1}\g(z)=g^{-1}\g(g)$.  Note
that $\g(gz^{-1})=gz^{-1}$ for all $\g\in\text{Gal}(F_s/F)$, so
$gz^{-1}\in G(F)$.  The result follows, since $gz^{-1}Nzg^{-1}=N'$.
\end{proof}

\begin{remark} If we assume that $\char\ k=0$, then the statement holds
for arbitrary semisimple $N$.  In this case, the field $F$ is perfect;
moreover, it has cohomological dimension $1$ because it is complete with
respect to a discrete valuation and has algebraically closed residue field
\cite[p. 97]{S2}.  This implies that $H^1(\text{Gal}(\overline F/F),L)=0$
for any connected linear group over $F$ \cite[I.9]{St}.  The theorem of
Springer and Steinberg again implies that the centralizer $Z$ is connected,
and the rest of the proof goes through as above.  \end{remark}

	A semisimple (resp. regular semisimple) element of $\frak{sl}(V_F)$
is just a (semisimple (resp. regular semisimple) element of
$\frak{gl}(V_F)$ with trace zero.  Also, two elements of
$\frak{sl}(V_{\overline F})$ are $SL(V_{\overline F})$-conjugate if and
only if they are $GL(V_{\overline F})$-conjugate.  Proposition 3.2 now shows
that the regular semisimple conjugacy classes of $\frak{sl}(V_F)$ under the
actions of $SL(V_F)$ and $GL(V_F)$ coincide.  We thus obtain as an
immediate corollary of theorem 2.4 and corollary 2.5 the following
classification of regular semisimple conjugacy classes $\eusm C$ in $\frak
{sl}(V_F)$ coming from almost commuting pairs $(N,f)$ with $f$ a
$Q'_{SL(V)}$-action:

\begin{theorem}  Let $\eusm C$ be a
regular semisimple (resp. nil regular semisimple) conjugacy class in
$\frak{sl}(V_F)$.  Then there is a representative $N\in \eusm C$ almost
commuting with a $Q'_{SL(V)}$-action if and only if the characteristic polynomial
$\text{char}_{\eusm C}(\mu)$ is a homogeneous polynomial of degree $n$ in
$\mu$ and $\pi^{q}$ for some nonzero (resp. positive) rational number $q\in
\frac1{n-1}\Z\cup \frac1n\Z$.  
\end{theorem}
\noindent Note that $c_1=\text{tr}\ \eusm C=0$ in $\char_{\eusm
C}(\mu)=\sum_{i=0}^n c_i \pi^{qi}\mu^{n-i}$, so the standard representative
for a conjugacy class given in the proof of theorem 2.4 lies in $\frak
{sl}(V_F)$.

	We are primarily interested in defining $k^*$-actions on the fixed
point varieties $\hat{\eusm P}_N$ when they are algebraic varieties. Recall
that a torus $T$ over $F$ is called {\it elliptic} if no nontrivial
cocharacter of $T$ is defined over $F$.  We say that a regular semisimple
element $N\in\frak{sl}(V_F)$ is elliptic if its centralizer
$Z_{SL(V_F)}(N)$ in $SL(V_F)$ is an elliptic torus.  Suppose that
$N\in\frak{sl}(V_F)$ is nil.  By a theorem of Kazhdan and Lusztig \cite{KL}, 
$\dim \hat{\eusm P}_N<\infty$ if and only if $N$ is regular semisimple.
Furthermore, $\hat{\eusm P}_N$ is an algebraic variety, i.e. has a finite
number of irreducible components precisely for $N$ elliptic.  (For positive
characteristic, we  note that the fixed point subvarieties 
considered below are algebraic varieties for an arbitrary algebraically
closed field; this follows from \cite{LS}.)

	It is known that a regular semisimple element $N$ of
$\frak{sl}(V_F)$ is elliptic if and only if $\char_N(\mu)\in F[\mu]$ is an
irreducible, separable polynomial \cite[proposition 3.10]{S}.  If $F$ has
characteristic zero, then for each positive
integer $m$, the field $F$ has a unique  extension of degree $m$ in
$\overline F$; it is generated by an $m$th root of $\pi$.  If we call this
cyclic extension $F(\pi^{1/m})$, we have $\overline F=\bigcup_{m=1}^\infty
F(\pi^{1/m})$.  In characteristic $p$, the structure of
$\overline F$ and the separable closure $F_s$ is much more complicated.
However, the previous theorem guarantees that any
elliptic, regular semisimple $N$ almost commuting with a
$Q'_{SL(V)}$-action has eigenvalues lying in
$\bigcup_{p\nmid m}
F(\pi^{1/m})\subset F_s$.  Therefore, if $N$ is nil-elliptic, the splitting field of
$\char_N$ is $F(\pi^{1/n})$ and $n$ can not be a multiple of the
characteristic.  By corollary 2.3, each eigenvalue is
homogeneous of valuation $q=s/n>0$ where $s$ is relatively prime to $n$.
Since $\pi^{qi}\notin F$ for $i\in[1,n-1]$, we have
$\char_N(\mu)=\mu^n-b\pi^s$.  Let $a\in k$ be some $n$th root of $b$.  Then
the roots of $\char_N$ over $\overline F$ are $\{a\zeta\pi^q:\zeta\
\text{an $n$th root of unity in $k$}\}$.  These are distinct if and only if
$a\ne 0$ and the characteristic of $k$ is either $0$ or relatively prime to
$n$.  This can be rewritten more simply as $(\char'\ k,n)=1$, where
$\char'\ k$ denotes the characteristic exponent of $k$.  (Recall that the characteristic exponent of a field is $1$ if the
field has characteristic $0$ and equals the characteristic otherwise
(\cite[5.1.5]{B2}).   Summing up, we have:

\begin{proposition} The nil-elliptic conjugacy
classes in $\frak{sl}(V_F)$ obtained from almost commuting pairs $(N,f)$
with $f$ a $Q'_{SL(V)}$-action are precisely those with $\char_{\eusm
C}(\mu)=\mu^n-b\pi^s$ where $b\in k^*$, $s>0$, $(s,n)=1$, and $(\char'\ k,n)=1$.  Fixing a
basis $\underline e$ for $V$, we have the standard matrix representative
\[N_{bs}=\begin{pmatrix} 0&1&0&0&\hdots&0\\ 0&0&1&0&\hdots&0\\
\vdots&\vdots&\ddots&\ddots&&\vdots\\ \vdots&\vdots&\ddots&\ddots&0&1\\ b\pi^s&0&0&\hdots&0&0
\end{pmatrix}. \]
Each $N_{bs}$ almost commutes with exponent $s$ with
the $Q'$-action $f$, diagonal with respect to $\underline e$, given
by $\nu(m,i)=nm-is$.  
\end{proposition}

Note that the field automorphism of $F$ given by
$\pi\mapsto b^{1/s}\pi$ (where $b^{1/s}$ is a fixed $s$th root of $b$)
induces an automorphism of $\frak{g}_F$ sending $N_{1s}$ to
$N_{bs}$. Accordingly, $\hat{\eusm{P}}_{N_{bs}}$ is isomorphic to
$\hat{\eusm{P}}_{N_{1s}}$, and we need only consider the case  $b=1$.  For
simplicity of notation, we set $N_s\overset{\text{def}}{=}N_{1s}$.

\begin{remark}  Lusztig and Smelt used different
representatives for these nil-elliptic conjugacy classes and almost commuting
actions \cite{LS}.  
\end{remark}

	For the remainder of the section, we assume that $(n,\char'\ k)=1$ and
fix a basis $\underline e$ of $V$.  By the Iwasawa decomposition for
$GL_n(F)$, any lattice $L$ has a basis $\underline z$ such that
$z_{\underline e}$ is
upper triangular.  Indeed, $L$ has a unique basis of the form
$z_j=\pi^{r_j}e_j+\sum_{i<j}z_{ij}e_i$ for $j\in[1,n]$, where $r_j\in\Z$
and $z_{ij}\in F$ is a terminating Laurent series with all nonzero terms of
degree strictly smaller than $r_i$.  The diagonal lattices (with
respect to $\underline e$) are the lattices $\L_{\bold
r}=A\pi^{r_1}e_1+\dots+A\pi^{r_n}e_n$ for $\bold r\in\bold Z^n$; they form a
countable closed subset $\eusm D$ of $\eusm L$.  Note that $\eusm
D_m\overset{\text{def}}{=}\eusm D\cap\eusm L_m=\{\L_{\bold r}\mid \sum_{i=1}^n
r_i=m\}$.  We call the set of lattices chains in $X_I$ composed of
diagonal lattices $\eusm D_I$.

	It is obvious that any $Q'$-action diagonal with respect to
$\underline e$ fixes every lattice in $\eusm D$.  The converse is not
true in general, but it does hold for the $Q'$-actions defined in the above
theorem.

\begin{lemma} Let $f$ be the $Q'$-action almost commuting with $N_{s}$ in
theorem 3.4.  The fixed point set of $f$ on $\eusm L$ is $\eusm D$.
\end{lemma}

\begin{proof}  Let $L$ be a lattice fixed by $f$, and choose a standard basis
$\underline z$ for $L$ as described above.  We will show that the
off-diagonal components of $z_j$ are zero by induction on $j$.  The claim is
vacuous for $j=1$.  Suppose that $j>1$.  Since $f(\l)z_j\in L$, we have
\begin{equation*}f(\l)z_j=\sum_{i'=1}^n a_{i'}z_{i'}\tag{3-1}\end{equation*} for some elements
$a_{i'}\in A$.  It is immediate that $a_{i'}=0$ for $i'>j$; if $i>j$ is the largest
index such that $a_i\ne 0$, then equating the coefficients of $e_i$ in this
expression gives $a_i=0$, a contradiction.  By comparing the $e_j$
terms, we now get $\l^{\nu(r_j,j)}\pi^{r_j}=a_j\pi^{r_j}$, i.e.
$a_j=\l^{\nu(r_j,j)}$.  If $i<j$, then the $e_i$ component of $z_{i'}$
vanishes unless $i'=i$ or $i'=j$ by inductive hypothesis.  Writing
$z_{ij}=\sum_{m\gg -\infty}^{r_i-1}z_{ij}^m\pi^m$ and equating the $e_i$
terms on each side of (3-1), we have

\[\sum_{m\gg
-\infty}^{r_i-1}\l^{\nu(m,i)}z_{ij}^m\pi^m=a_i\pi^{r_i}+\l^{\nu(r_j,j)}\sum_{m\gg
-\infty}^{r_i-1}z_{ij}^m\pi^m.\]
 Therefore, \[a_i=\pi^{-r_i}\sum_{m\gg
-\infty}^{r_i-1}(\l^{\nu(m,i)}-\l^{\nu(r_j,j)})z_{ij}^m\pi^m\in A.\]  Since
$m-r_i<0$, $(\l^{\nu(m,i)}-\l^{\nu(r_j,j)})z_{ij}^m=0$ for each $m$.  If
$\l^{\nu(m,i)}-\l^{\nu(r_j,j)}=0$, then $n(r_j-m)=(j-i)s$, a contradiction
since $(s,n)=1$ and $n$ does not divide $j-i\in[1,n-1]$.  Thus, $z_{ij}^m=0$
for all $m$, i.e. $z_j=\pi^{r_j}e_j$.  
\end{proof}

\begin{remarks} \begin{mrk} \item Using a similar argument,
it is possible to show that if $f$ is a diagonal $Q'$-action almost
commuting with a nonnilpotent element $N\in\frak{gl}_n(F)$ and if $(N,f)$ has
loop index $s/n$ with $(s,n)=1$, then the fixed point set of $f$ on $\eusm
L$ is $\eusm D$ \cite[proposition 3.12]{S}.  \item It is an immediate consequence that $\eusm D$ is closed in
$\eusm L$.  \end{mrk}
\end{remarks}

	  For any nil-elliptic $N\in\frak{sl}(V_F)$, $\eusm
D\cap\eusm L_{m,N}$ is an algebraic variety and hence has a finite number of
components.  Since 
$\eusm D$ is discrete, the proposition then shows that the fixed point sets
$(\hat{\eusm P}_{N_{s}})^{k^*}$ for the $k^*$-action $f$ of proposition 3.4
are all finite.  A theorem of Bialynicki-Birula \cite{B-B1} implies that $\chi(\hat{\eusm P}_{N_{s}})=\chi((\hat{\eusm
P}_{N_{s}})^{k^*})$, so the calculation
of these Euler characteristics has been reduced to a problem in
combinatorics.

	We can decompose the computation into two complementary parts.  Take
$I\subsetneq [0,n-1]$ and let $J=[0,n-1]\setminus\{m\}$ be a maximal proper subset of
$[0,n-1]$ containing $I$.  The fiber bundle $\hat{\eusm P}_I\overset{\rho}{\to}\hat{\eusm P}_J$ with fiber $\eusm P$ induces a map $\hat{\eusm
P}_{I,N}\to\hat{\eusm P}_{J,N}$.  Upon passing to points fixed by $f$, we see that
\begin{equation*}\chi(\hat{\eusm P}_{I,N})=\chi((\hat{\eusm P}_{I,N})^{k^*})=\sum_{z\in (\hat{\eusm
P}_{J,N})^{k^*}} |(\rho^{-1}(z))^{k^*}|.\tag{3-2}\end{equation*} 
Suppose $z\in (\hat{\eusm P}_{J,N})^{k^*}$ corresponds to the
lattice $L\in\eusm L_m$.  Let $\bar{N}\in\frak{sl}(L/\pi L)$ be the
nilpotent endomorphism induced by $N$.  The summand may then be interpreted
as the Euler characteristic of $\eusm P_{\bar{N}}$, a fixed point variety on
a classical flag
manifold.  We can thus approach the problem by examining
the fixed points in the base and fibers of $\rho$ separately.

	We will first determine the fixed point sets $(\eusm
L_{m,N_{s}})^{k^*}$.  Let $\L$ be a lattice of valuation $m$ fixed by $f$
and $N_{s}$.  Then $\L$ must equal $\L_{\r}$ for some
$\r\in\Z^n$ with $\sum_{i=1}^n r_i=m$.  The inclusion $N_{s}\L\subset\L$ implies that
$\pi^{r_i}e_{i-1}\in A\pi^{r_{i-1}}e_{i-1}$ for $i\in[2,n]$ and
$\pi^{r_1+s}\in A\pi^n$.  Thus, $(\eusm L_{m,N_{s}})^{k^*}=\{\L_{\r}\mid \r\in
R_{sm}\}$ where 
\begin{equation*}R_{sm}=\{\r\in\Z^n\mid \sum_{i=1}^n r_i=m \ \text{and $r_1\le
r_2\le\dots\le r_n\le r_1+s$}\}.\tag{3-3}\end{equation*}

	To determine the size of $R_{sm}$, we will use an argument due to
Lusztig and Smelt \cite{LS}.  For $i\in[1,n]$, set $\z_i=\d_{i1}s$ and
define a function $\Th_s:\Z\to\Z$ by
$\Th_s(l)=\sum_{i=1}^n(l+i)(r_i-r_{i-1}+\z_i)$ (where we view $\r$ as a
periodic sequence indexed by $\Z$).  We have \begin{equation*}\begin{split}
\Th_s(l)=\Th_s(0)+l\sum_{i=0}^{n-1}(r_{i+1}-r_i+\z_i)&=nr_n-\sum_{i=1}^{n}r_i+s+ls\\
&=-m+(l+1)s+nr_n\\ &\equiv -m+(l+1)s\pmod n.\end{split}\tag{3-4} \end{equation*} Hence,
$\Th_s$ induces a map $\Z/n\Z\to\Z/n\Z$, and it is in fact a bijection
because $s$ and $n$ are relatively prime.  Set $C^n_s=\{\c\in\N^n\mid
\sum_{i=1}^nc_i=s\}$ and $C^n_{sm}=\{\c\in C^n_s\mid
\sum_{i=1}^n(l+i)c_i\equiv -m+(l+1)s\pmod n\}$.  Usually the superscript
will be suppressed from the notation.  Just as in (3-4), the
congruence will hold for each integer if it holds for a single $l$.
The group $\Z/n\Z$ acts on $C_s$ via cyclic permutations of the coordinates; let
$\tilde{C}_s$ denote the set of orbits.

\begin{lemma} There are bijections
$R_{sm}\overset{\phi}{\to} C_{sm}\overset{\psi}{\to}\tilde{C}_s$.
\end{lemma}

\begin{proof} Set $\phi(\r)_i=r_i-r_{i-1}+\z_i$.  If $\phi(\r)=\phi(\r')$,
then equality of the last $n-1$ components implies that $\r'=\r+t\bold 1$
where $t\in\Z$ and $\bold 1$ is the constant $n$-tuple $(1,\dots,1)$.  The
sums of the coefficients are equal, so $t=0$.  Thus, $\phi$ is injective.
Now take $\c\in C_{sm}$.  Set \[r_n=\frac{m-s+\sum_{i=1}^nic_i}{n} \text{ and
$r_j=r_n-(c_{j+1}+\dots+c_n)\,$ for $j\in [1,n-1]$.}\]  These are integers
by definition of $C_{sm}$.  We have
\[\sum_{i=1}^nr_i=nr_n-\sum_{i=1}^n(i-1)c_i=m-s+\sum_{i=1}^nc_i=m.\] This
proves surjectivity.

	Let $\psi$ be the restriction of the orbit map to $C_{sm}$.
The $n$-cycle $\tau=(1,2,\dots, n)$ generates the $\Z/n\Z$ action.  Note that for any
$\c\in C_s$, $\sum_{i=1}^ni(\tau\cdot\c)_i=\sum_{i=1}^n(1+i)c_i-nc_n\equiv
s+\sum_{i=1}^n ic_i\pmod n$.  This means that within the orbit of $\c$, the sum runs
through the residue classes mod $n$, i.e. precisely one representative
lies in $C_{sm}$.  Therefore, $\psi$ is bijective.  
\end{proof}

\begin{corollary} For any $m\in\Z$, \[\chi(\eusm
L_{m,N_{s}})=\frac{(n+s-1)!}{n!s!}.\]  
\end{corollary}

\begin{proof} The proof of the lemma shows that $\Z/n\Z$ acts freely on
$C_s$.  Consequently, $\chi(\eusm
L_{m,N_{s}})=|R_{sm}|=|\tilde{C}_s|=|C_s|/n$.  It now suffices to
recall the well known fact that
$|C_s|=\binom{n+s-1}{n-1}=\frac{(n+s-1)!}{(n-1)!s!}$; it is just the
number of ways of putting $s$ indistinguishable balls into $n$ boxes.
\end{proof}

	We have just shown that $|\tilde{C}^n_s|=\frac{(n+s-1)!}{n!s!}$.
This expression is symmetrical in $n$ and $s$, so we have a bijection
$\tilde{C}_s\tilde{\to}\tilde{C}^s_n$.  We
can describe this bijection explicitly.  The orbit $\c\in\tilde{C}^n_s$ may be
interpreted as an arrangement of $s$ indistinguishable balls in $n$ boxes
placed around a circle.  Viewing a box as two boundary walls, we see that this
is just a configuration of $s$ balls and $n$ walls up to cyclic permutation
of the walls.  However, because $n$ and $s$ are relatively prime, this is equivalent to an
arrangement of the $n+s$ objects up to cyclic permutation.  Let a cell be
a ``dual box'', where the boundaries are now balls and the contents of the
cell are walls.  The analogous
argument now gives a way of placing $n$ walls in $s$ cells around a circle
and hence defines the desired bijection.  More concretely, number the boxes
and walls clockwise around the circle so that the $(t-1)$st and $t$th walls
bound the $t$th box.  Label the balls and cells similarly, starting with
the first ball counterclockwise of wall $1$.  Letting $B_k$ denote the set
of walls contained in the $k$th cell, we have $b_k=|B_k|$.  The cell
numbering guarantees that for walls $t,t'$ within a given cell, $t<t'$ if
and only if the path from $t$ to $t'$ staying in the cell runs clockwise.

	We now turn to lattice chains of type $J'$ where $J'$ is a nonempty
subset of $[0,n-1]$.  Let $l$ be the number of elements in $J'$.  Choose
$m\in J'$ and let $J=\{j_1< \dots<j_l=n\}\subseteq
[1,n]$ be the unique set such that
$J'=m+((J\setminus\{n\})\cup\{j_0=0\})\pmod n$.  Set $p_i=j_i-j_{i-1}$.  A
chain in $X_{J'}$ is uniquely determined by its component lattices in
$\eusm L_{m+j}$ for $j\in J$.  Thus, $(X_{J',N_{s}})^{k^*}=\{(\L_{\r^j})_{j\in
J}\mid (\r^0,\dots,\r^{j_l})\in Y_{Jsm}\}$ where
\[Y_{Jsm}=\{(\r^0,\dots,\r^{j_l})\mid \r^{j_i}\in R_{s,m+j_i} \text{ and }
\r^0\le\r^{j_1}\le\dots\le\r^{j_{l-1}}\le\r^n=\r^0+\bold 1\}.\]  The partial
order on vectors is defined componentwise.

	Fix $\r\in R_{sm}$ and let $Y_{Jsm}(\r)$ denote those chains in
$Y_{Jsm}$ with $\r^0=\r$.  Given any chain $\hat{\r}$ in this fiber, $r_t\le
r^{j_1}_t\le\dots\le r^{j_l}_t=r_t+1$ for each $t\in [1,n]$.  Hence, there
exists a unique $\s(t)=j_i\in J$ such that $r^{j_{i-1}}_t=r_t$ and
$r^{j_i}_t=r_t+1$.  The map $\s:[1,n]\to J$ partitions $[1,n]$ into subsets
$\s^{-1}(j_i)$ of size $p_i$ and uniquely determines $\hat{\r}$.  Let $\b\in
C^s_n$ be a representative of the element of $\tilde{C}^s_n$ corresponding to
$\r$ through the bijections $R_{sm}\to\tilde{C}_s\to\tilde{C}^s_n$.  We
define a set of functions 
\[ 
Q'_{Js}(\b)=\{\s:[1,n]\to J\bigm|  |\s^{-1}(j)|=p_j\text{ and $\s$ is
nonincreasing on each $B_k$}\}.\] 
Note that although this set depends on a
choice of labeling of the cells $B_k$, any other labeling (not just those
obtained through cyclic permutations) or even any decomposition of $[1,n]$
into $s$ disjoint subsets of sizes $b_k$ gives a set of the same size.

\begin{proposition} The map $\hat{\r}\mapsto\s$ defines a bijection
$Y_{Jsm}(\r)\tilde{\to} Q'_{Js}(\b)$.  
\end{proposition}

\begin{proof} Take $\hat{\r}\in Y_{Jsm}$.  Note
that $t-1$ and $t$ are in $B_k$ if and only if $c_t=0$; $c_t\ne 0$ means that box $t$
contains a ball, separating walls $t-1$ and $t$ into different cells.
Suppose that $t-1$ and $t$ are in $B_k$, but $\s(t)>\s(t-1)=j_i$.  Then $r^{j_i}_t=r_t$
and $r^{j_i}_{t-1}=r_{t-1}+1$, implying that
$c^{j_i}_t=r_t-(r_{t-1}+1)=c_t-1=-1$, a contradiction.  Moreover, for each
$j\in J$, $\sum_{t=1}^n r^j_t=m+j$, so $r^{j_{i-1}}_t<r^j_t$ for exactly
$p_i$ elements of $[1,n]$.  Thus, the map sends $Y_{Jsm}$ into
$Q'_{Js}(\b)$.

The map is clearly injective.  Take $\s\in Q'_{Js}(\b)$ and define $r^{j_i}$ recursively
by 
\[r^{j_i}_t=\begin{cases} r^{j_{i-1}}_t+1&\text{if $t\in\s^{-1}(j_i)$}\\
r^{j_{i-1}}_t&\text{if $t\notin\s^{-1}(j_i)$}.\end{cases}\] 
We will show that
$\r^{j_i}\in R_{s,m+j_i}$ by induction on $i$.  This is trivial for
$i=0$.  Now suppose that $i>0$ and that the claim is true for indices
smaller than $i$.  In particular, $\c^{j_{i'}}=\phi(\r^{j_{i'}})$ is
well-defined for all $i'<i$.  The coordinates of the $n$-tuple $\r^{j_i}$
sum to $m+j_i$ because $\sum_{t=1}^n
r^{j_i}_t=\sum_{t=1}^n
r^{j_{i-1}}_t+|\s^{-1}(j_i)|=m+j_{i-1}+p_i=j_i$.  Also, observe
that if $r^{j_i}_t+\z_t<r^{j_i}_{t-1}$, then
$r^{j_{i-1}}_t+\z_t=r^{j_{i-1}}_{t-1}$, i.e. $c^{j_{i-1}}_t=0$.  Thus,
$\s(t-1)=j_i$, and since
$c^{j_{i'}}_t$ only decreases where $r^{j_{i'}}_{t-1}$ increases, namely
at the $i$th step, $c^{j_{i'}}_t=0$ for each $i'<i$.  In
particular, $c_t=0$.  But this means that $t$ and $t-1$ are in the same
cell, so $j_i=\s(t-1)\ge\s(t)$.  This is a contradiction; if
$\s(t-1)=\s(t)$, then $0=r_t^{j_{i-1}}-r_{t-1}^{j_{i-1}}+\zeta_t=r_t^{j_{i}}-r_{t-1}^{j_{i}}+\zeta_t<0$ and if
$\s(t-1)>\s(t)=j_{i'}$, then $c^{j_{i'}}_t=1$ for some $i'<i$.  This
proves surjectivity.  
\end{proof}

\begin{definition} Let $\bold d=(d_0,d_1,\dots,d_l)$ be a
nondecreasing sequence and $t$ a positive integer.  We call a matrix
$(q_{ik})\in M_{lt}(\N)$ a $(\bold d,t)$-{\it intersection matrix} if
$\sum_{k=1}^t q_{ik}=d_i-d_{i-1}$ for $i\in [1,l]$.  We denote the set of
$(\bold d,t)$-intersection matrices by $Q_{\bold
d t}$.  If $\bold d$ is
strictly increasing, then the $d_i-d_{i-1}$'s are determined by
$D=\{d_0,\dots,d_l\}$.  In this case, $Q_{Dt}\overset{\text{def}}{=}Q_{\bold
d t}$ is called the set of $(D,t)$-intersection matrices.  For any $\b\in
C^t_{d_l-d_0}$, let $Q_{\bold d t}(\b)\subseteq Q_{\bold d t}$ be the set
of $(\bold d,t)$-intersection matrices with $\sum_{i=1}^l q_{ik}=b_k$ for each
$k\in [1,t]$.  The group $\Z/t\Z$ acts freely on $Q_{\bold d t}$ by cyclic
permutation of the columns; we call the set of orbits $\tilde{Q}_{\bold d
t}$.  \end{definition}

\begin{lemma} \[|Q_{\bold d
t}|=\prod_{i=1}^l\binom{s+d_i-d_{i-1}-1}{d_i-d_{i-1}}.\] 
\end{lemma}

\begin{proof} The $i$th row of a $(\bold d,t)$-intersection matrix is just an
element of $C^t_{d_i-d_{i-1}}$, so we have a bijection $Q_{\bold d
t}\to\prod_{i=1}^l C^t_{d_i-d_{i-1}}$.  The result follows, since
$|C^t_{d_i-d_{i-1}}|=\binom{t+d_i-d_{i-1}-1}{d_i-d_{i-1}}$ by the proof
of corollary 3.7.  
\end{proof}

	We will describe the Euler characteristics of $X_{J,N_{s}}$ in
terms of $(J,s)$-in\-ter\-sec\-tion matrices.  Define a map $Q'_{Js}(\b)\to
Q_{Js}(\b)$ by $\s\mapsto (|\s^{-1}(j_i)\cap B_k|)$.  (This is the
motivation for the terminology.)

\begin{proposition} The map $Q'_{Js}(\b)\tilde{\to} Q_{Js}(\b)$ is a
bijection.  
\end{proposition}

\begin{proof} Fix an intersection matrix $Q=(q_{ik})\in Q_{Js}(\b)$, and suppose
$\s$ is mapped to $Q$.  The function $\s$ is nonincreasing on $B_k$, so the top $q_{1k}$
elements of $B_k$ are in $\s^{-1}(1)$, the next highest $q_{2k}$ are in
$\s^{-1}(2)$, and so on.  Thus, $\s$ is completely determined by its
intersection matrix.  Conversely, this procedure clearly defines $\s\in
Q'_{Js}(\b)$ mapping to $Q$.  
\end{proof}

	We can apply this proposition to obtain a combinatorial description
of the Euler characteristic of $\eusm P_N$ for any classical partial flag
variety $\eusm P$ and any nilpotent $N\in\frak{sl}_n(k)$.  Nilpotent
conjugacy classes in $\frak{sl}_n(k)$ are determined by the Jordan blocks
of $N$, in other words, by a partition of $n$, say $n=n_1+\dots+n_t$ where
each $n_i$ is positive.

\begin{theorem}   Assume that $(\char'\ k,n)=1$.  Let $N\in\frak{sl}_n(k)$ be a nilpotent element in
the conjugacy class corresponding to the partition $n=n_1+\dots+n_t$.  Let
$s$ be any integer greater than $t$ and relatively prime to $n$.  Define
$\b\in C^s_n$ by $b_i=n_i$ for $i\in [1,t]$ and $0$ otherwise.  Then for each
$I\subseteq [1,n-1]$, $\chi(\eusm F_{I,N})=|Q_{I\cup\{n\},s}(\b)|$.  In
particular, $\chi(\eusm B_N)=\dfrac{n!}{n_1!\cdots n_t!}$. 
\end{theorem}

\begin{proof} Let $\r\in R_{s0}$ correspond to the orbit of $\b$.  Since the
fiber of the fixed point set over $\r$ corresponds bijectively to
$Y_{I\cup\{n\},s0}(\r)\tilde{\to} Q_{I\cup\{n\},s}(\b)$, it will
suffice to show that the nilpotent map $\bar{N}_{s}$ induced on
$W=\L_{\r}/\pi\L_{\r}$ has Jordan block structure given by the partition of
$n$.

	The images $w_i$ of $\pi^{r_i}e_i$ form a basis for $W$, and the
matrix for $\bar{N}_{s}$ with respect to this basis has the same form as
does $N_{s}$.  Thus, only the $(n,1)$
and $(i,i+1)$ coefficients can possibly be nonzero.  The image of $w_i$ under $\bar{N}_{s}$ is zero if and only if
$r_i-r_{i-1}+\z_i\ge 1$, i.e. exactly when $c_i>0$.  Suppose the $i$th and
$i'$th columns are nonzero with vanishing intervening columns (in
increasing order with $1$ following $n$).  Then $c_i$ and $c_{i'}$ are
positive, and there is a cell containing precisely the walls $i$ through $i'-1$.  But
a Jordan block of size $l$ is given by a string of $l-1$ consecutive
nonzero columns surrounded by two zero columns, implying that Jordan blocks
correspond to cells.  The $b_k$'s are the cell sizes, so the claim is
proved.

	Since $\eusm B=\eusm F_{[1,n-1]}$, we have $\chi(\eusm
B_N)=|Q_{[1,n],s}(\b)|$.  Each row of a $([1,n],s)$-in\-ter\-sec\-tion matrix
sums to $1$, implying that the only nonzero entry in any row is a $1$.
Such a matrix is in the fiber over $\b$ if and only if there are exactly
$b_k$ $1$'s in the $k$th column.  Thus, an element of $Q_{[1,n],s}(\b)$ is
determined by arbitrarily selecting subsets of sizes $b_1,\dots,b_k$ from
the $n$ rows.  This can be done in $\frac{n!}{b_1!\cdots b_s!}=\frac{n!}{n_1!\cdots n_t!}$ ways.
\end{proof}

\begin{remark}  The smallest integer which can be used for
$s$ in the
above theorem for every nilpotent conjugacy class is $n+1$.
\end{remark}

	The bijection between the fiber $Y_{Jsm}(\r)$ and $Q_{Js}(\b)$
extends to a one-to-one correspondence between the entire fixed point set
and all orbits of $(J,s)$-intersection matrices.

\begin{proposition} There is a bijection
$Y_{Jsm}\overset{\a}{\to}\tilde{Q}_{Js}$.  
\end{proposition}

\begin{proof} We define the function $\a$ on $Y_{Jsm}(\r)$ to be the composition of the
maps $Y_{Jsm}(\r)\to Q'_{Js}(\b)\to Q_{Js}(\b)\to\tilde{Q}_{Js}$.  The
$s$-tuple $\b\in
C^s_n$ is only defined up to cyclic permutation, but this ambiguity
disappears upon passing to orbits.

	Suppose $\a(\hat{\r})=\a(\hat{\r}')$.  Note that $Q_{Js}(\b)$ and
$Q_{Js}(\b')$ have disjoint images in the orbit space unless $\b$ and $\b'$
are cyclic permutations of each other.  Hence, $\hat{\r}$ and $\hat{\r}'$ must
be in the same fiber $Q_{Js}(\b)$.  Propositions 3.8 and 3.10 now imply that
$\hat{\r}=\hat{\r}'$, so $\a$ is injective.

	Given $Q=(q_{ik})\in Q_{Js}$, define $\b$ by $b_k=\sum_{i=1}^l
q_{ik}$.  These are nonnegative integers, satisfying $\sum_{k=1}^s
b_k=\sum_{k=1}^s\sum_{i=1}^l q_{ik}=\sum_{i=1}^l\sum_{k=1}^s
q_{ik}=\sum_{i=1}^l p_i=n$. This shows that $\b\in C^s_n$ and $Q\in
Q_{Js}(\b)$.  It is now clear that the corresponding orbit is in the image
of $\a$, proving surjectivity.  
\end{proof}

	We obtain the following theorem on Euler characteristics of fixed
point varieties on affine flag manifolds:

\begin{theorem} Assume that $(n,\char'\ k)=1$, and let $s$ be a positive
integer relatively prime to $n$.  Let
$I=\{j_1<\dots<j_l\}\subseteq [0,n-1]$ be a nonempty set of size $l$ with
$p_i=j_{i+1}-j_i$ (where $j_{l+1}=j_1+n$).  Then
\[\chi(X_{I,N_{s}})=s^{-1}\prod_{i=1}^l\binom{s+p_i-1}{p_i}.\] In
particular, $\chi(\hat{\eusm B}_{N_{s}})=s^{n-1}$ and if $\eusm K$ is any
variety of maximal parahoric subalgebras, $\chi(\hat{\eusm
K}_{N_{s}})=\dfrac{(n+s-1)!}{n!s!}$.  
\end{theorem}

\begin{proof} Let $J\subseteq [1,n]$ be the unique set such that
$I=j_1+((J\setminus\{n\})\cup\{0\})\pmod n$.  Note that
$J=\{j_2-j_1.\dots,j_l-j_1,n\}$, so the definition of the $p_i$'s given in
the statement coincides with the earlier definition.  Since there is a
bijection $Y_{Jsj_1}\tilde{\to}\tilde{Q}_{Js}$ by proposition 3.12, we have
\[\chi(X_{I,N_{s}})=\chi((X_{I,N_{s}})^{k^*})=|Y_{Jsj_1}|=|\tilde{Q}_{Js}|=\frac{|Q_{Js}|}s.\] The formula now follows from lemma 3.9.

	The varieties of maximal parahorics correspond to lattice chains of
type $I$ with $I$ a singleton.  In this case, $l=1$ and $p_1=n$, and the
Euler characteristic is
\[\frac{(n+s-1)}{n!(s-1)!s}=\frac{(n+s-1)!}{n!s!}.\] We thus recover
corollary 3.7.  The full affine flag variety $\hat{\eusm B}$ is isomorphic to
$X_{[0,n-1]}$.  This gives $l=n$ and $p_i=1$ for each $i$, so $\chi(\hat{\eusm
B}_{N_{s}})=s^{-1}\dbinom{s}{1}^n=s^{n-1}$.  \end{proof}

\begin{remarks} \begin{mrk} \item  Let $\lambda_1,\dots,\lambda_{n-1}$ be
the roots of $\frak{sl}_n(k)$ and $\lambda_0$ the affine root, and consider
the subset of the positive Weyl chamber in an apartment of the affine building given by the conditions
$\lambda_i\ge 0$ for $i\in[1,n-1]$ and $\lambda_0\le s$.  The Euler
characteristic of $\hat{\eusm{P}}_{N_{s}}$ is then  just the number of simplices of type corresponding
to the parahoric class in this region.
\item For $s=1$, the theorem implies that the
Euler characteristic is $1$ for each parahoric while if $n$ is odd and
$s=2$, the Euler characteristics are given by \[\chi(X_{I,N_{b2}})=\frac
12\prod_{i=1}^l(p_i+1).\]  \item Lusztig and
Smelt obtained the formulas for $\chi(\hat{\eusm B}_{N_{1s}})$ and
$\chi(\eusm L_{0,N_{1s}})$ \cite{LS}.  \end{mrk} \end{remarks}

Finally, we recall a theorem due to Lusztig and Smelt:
\begin{theorem}  The varieties $\hat{\eusm{P}}_{N_{s}}$ have 
stratifications with finitely many strata, each of which is isomorphic to
an affine space over $k$.
\end{theorem}
\begin{proof}  Since $N_{s}$ is conjugate to the nil-elliptic
 element $T_s$ considered by Lusztig and Smelt, their result carries over
 immediately\cite{LS}.  (Their argument, given for Iwahori and maximal parahoric
 subalgebras applies to the other parahoric subalgebras
 as well.)
\end{proof}

\begin{remark}  The particular $k^*$-action considered by Lusztig and Smelt
 came from the $Q'$-action $f(\lambda)\pi^me_i=\lambda^{mn-i}\pi^me_i$.
 The cells of the Bruhat decomposition for $\hat{\eusm{P}}$ for this choice
 of ordered basis are just the
 cells of the Bialynicki-Birula plus decomposition for $f$.  For example,
 consider the Bruhat cells of $\eusm{L}_m$.  Set
 $i,j\in[1,n]$ with $i\ne j$, set $<i,j>=0$ if $i<j$ and $1$ otherwise.
 The cells are indexed by
 $\mathbf{r}\in\mathbf{Z}^n$ such that $\sum_{i=1}^n r_i=m$, and they
 consist of those lattices with an $A$-basis $z_1,\dots,z_n$ of the form
 $z_j=\pi^{r_j}(e_j+\sum_{i\ne j}(\sum_{l\ge 0}
 a^l_{ij}\pi^{l+<i,j>})e_i)$, where $a^l_{ij}\in k$ for
 $i\ne j, l\ge 0$ and $a^l_{ij}=0$ if $l\ge r_i-r_j-<i,j>$.  It is easy to
 see that $f$ preserves the Bruhat cells and acts on the coordinates
 $a^l_{ij}$ for a given cell by sending them to
 $\lambda^{nl+n<i,j>+j-i}a^l_{ij}$.  The coefficient of $\lambda$ is
 always positive, so the limit as $\lambda$ goes to zero of $f(\lambda)L$ for any lattice
 in this cell is the diagonal lattice with $A$-basis
 $\pi^{r_1}e_1,\dots\pi^{r_n}e_n$.  Since the cells of Lusztig and Smelt's
 stratification of $\hat{\eusm{P}}_{T_s}$ are just the intersections of the
 fixed point variety with the Bruhat cells, we see that the stratification
 coincides with the plus decomposition.
\end{remark}

\begin{corollary} The \'{e}tale cohomology of the varieties
 $\hat{\eusm{P}}_{N_{s}}$ vanishes in odd degrees.

\end{corollary}

\section{$\lowercase{k}^*$-actions on fixed point varieties of affine flag
manifolds: The symplectic algebras}

	We now consider the simple root systems $C_n$.  Let $V$ be a
$2n$-dimensional vector space over $k$ endowed with a nondegenerate,
alternating, bilinear
form $\<\ ,\ \>$.  The symplectic group $Sp(V)$ is the
group of linear automorphisms of $V$ preserving the form.  The
corresponding Lie algebra is denoted $\frak{sp}(V)$.  A basis $\underline
e=(e_1,\dots,e_{2n})$ of $V$ is called {\it symplectic} if the subspaces
$ke_1\oplus\dots\oplus ke_n$ and $ke_{n+1}\oplus\dots\oplus ke_{2n}$ are
isotropic, i.e. the restriction of the form to these subspaces is
identically $0$, and
$\<e_i,e_{n+j}\>=\delta_{ij}$ for $i,j\in [1,n]$.  It is well known that $V$
admits such a basis.  Define a matrix
\[J=\begin{pmatrix} 0 & I_n\\ -I_n & 0 \end{pmatrix}.\]
A choice of symplectic basis identifies $Sp(V)$ with $Sp_{2n}(k)=\{g\in
GL_{2n}(k)\mid g^tJg=J\}$, the semisimple, simply connected,  algebraic group of type $C_n$.  If $n=1$, this is
just $SL_2(k)$.  The symplectic Lie algebra with respect to this
basis $\frak{sp}_{2n}(V)$ is the space of $2n\times 2n$ matrices $M$
satisfying $M^tJ+JM=0$; equivalently, 
\begin{equation*} M=\begin{pmatrix} M_{11} & M_{12}\\ M_{21}
& M_{22} \end{pmatrix}\tag{4-1}\end{equation*} 
where $M_{22}=-M_{11}^t$ and the
off-diagonal $n\times n$ blocks are symmetric.

	The symplectic basis $\underline{e}$ determines a standard
Borel subalgebra $\frak b$ consisting of the matrices $M\in \frak{sp}_{2n}(k)$
such that $M_{21}=0$ and $M_{11}$ is upper triangular, using the notation
of (4-1). (The subalgebra of upper triangular matrices in $\sp(k)$ is not a
Borel subalgebra.)  Let $\hat{\frak b}\subset\sp(A)$ be the
corresponding standard Iwahori subalgebra.  We define the algebras
$\hat{\frak q}_i\subset\sp(F)$ for $i\in[0,n]$ by
\begin{align*}
 &\qquad\begin{matrix}  \negthickspace i&n-i&\!\quad i&\thickspace n-i\end{matrix}\\ 
\hat{\frak q}_i=\sp(F)\ \cap &\begin{pmatrix}
A&A&\pi^{-1}A&A\\
\pi A&A&A&A\\
\pi A&\pi A&A&\pi A\\
\pi A&A&A&A \end{pmatrix} \begin{matrix} i\\n-i\\i\\n-i
\end{matrix} 
\end{align*}
These are the maximal parahoric subalgebras $\hat{\frak q}_i$ containing our fixed Iwahori subalgebra.  Their
stabilizers $\hat{Q}_i\subseteq Sp_{2n}(F)$ have the same block form.  The parahoric subalgebras
containing $\hat{\frak b}$ are now given by $\hat{\frak{p}}_I=\cap_{i\notin I}\hat{\frak
q}_i$ with stabilizer $\hat{P}_I=\cap_{i\notin I}\hat{Q}_i$ for each
proper subset $I$ of $[0,n]$; this notation is compatible with notation of
section 2 for a suitable choice of simple reflections for $\hat{W}$, indexed
by $[0,n]$.  The affine flag manifolds $\hat{\eusm P}_I$ for
$Sp(V_F)$ are identified with the quotient spaces $Sp(V_F)/\hat{P}_I$.

	We will again interpret the spaces of parahoric subalgebras as
chains of $A$-lattices in $V_F$.  Let $L\subset V_F$ be a lattice.  The set $L^*=\{v\in V_F\mid
\<v,L\>\subseteq A\}$ is called the {\it dual lattice} to $L$; it is a
lattice with valuation $-v(L)$.

\begin{definition} A lattice $L$ is called {\it symplectic}
if the set of lattices generated from $L$ by homothety and duality (or
equivalently, the set of lattices homothetic to $L$ or $L^*$) is totally
ordered by inclusion, i.e. a chain.  The variety of symplectic lattices is
denoted $\Ls$.
\end{definition}

	The symplectic group $Sp(V_F)$ acts on $\Ls$, and the orbits of
this action are the subvarieties $\Ls_m\overset{\text{def}}{=}\Ls\cap\eusm
L_m$.  We can apply the theory of alternating forms on free modules over a
principal ideal domain (see theorem IX.5.1 in \cite{B1}) to obtain a more
explicit characterization of these orbits.  For $m\ge 0$, write $m=qn+r$
with $r\in [0,n-1]$.  The variety $\Ls_m$ is the set of lattices contained
in their duals and admitting a basis $\underline{z}$ such that
$\bigoplus_{i=1}^n Az_i$ and $\bigoplus_{i=n+1}^{2n}Az_i$ are isotropic
submodules and $\<z_i,z_{n+j}\>$ equals $\delta_{ij}\pi^{q}$ for
$i\in[1,n-r]$ and $\delta_{ij}\pi^{q+1}$ for $i> n-r$.  For $m\le 0$,
$\Ls_m$ is the set of lattices whose duals are in $\Ls_{|m|}$.  In
particular, $\Ls_0$ is the space of self-dual lattices.

	As in the case of $SL(V_F)$, if $m\equiv m'\pmod{2n}$, then
$\Ls_m\cong\Ls_{m'}$ as an $Sp(V_F)$-variety.  In addition, the
correspondence $L\mapsto L^*$ defines an $Sp(V_F)$-isomor\-phism
$\Ls_m\to\Ls_{-m}$.  Thus, by identifying a symplectic lattice $L$ with the
totally ordered set $\{\pi^qL,\pi^{q'}L^*\mid q,q'\in\Z\}$, we can view
$\Ls_m$ as a variety of chains of symplectic lattices.  More generally, we
make the following definition:

\begin{definition} Let $I$ be a nonempty subset of $[0,n]$.  The space
of (partial) symplectic lattice chains of type $I$ is the set
\begin{multline*}
\Xs_I=\{(L_i)\mid i\in I\,\cup -I\negmedspace\pmod{2n}, L_i\in\Ls_i,\\ L_i\supset L_j\iff
i<j,\,
L_{i+2n}=\pi L_i,\ \text{and $L_{-i}=L^*_i$}\}.
\end{multline*} 
\end{definition}
\noindent The $Sp(V_F)$-action on the $\Ls_m$'s extends to an action on the
$\Xs_I$'s.  Note that $\Xs_I$ is a subvariety of
$X_{I\cup(2n-I)\setminus\{2n\}}$, the set of lattice chains of type
$I\cup(2n-I)\setminus\{2n\}$.

	An argument similar to that sketched in section 3 shows that the
varieties $\Xs_I$ are the affine flag manifolds for $Sp(V_F)$ (see
proposition 4.6 and theorem 4.8 in \cite{S}): 

\begin{theorem} For each
proper subset of $[0,n]$, $\hat{\eusm P}_I\cong\Xs_{I^c}$ as an
$Sp(V_F)$-variety.  Moreover, if $0\in I$, then the
fiber over the self-dual lattice $L$ of the projection
$\Xs_I\to\Xs_{\{0\}}=\Ls_0$ is isomorphic as an $Sp(L/\pi L)$-variety to
$\eusm F^{Sp}_{I\setminus\{0\}}(L/\pi L)$, the variety of isotropic flags
of type $I\setminus\{0\}$ in $L/\pi L$.
\end{theorem}

	 Let $\tilde{Q}'_{Sp}=GSp(V)\times k^*$, where $GSp(V)$ is the
group of symplectic similitudes of $V$.  We call a $Q'$-action $f$ a
$\tilde{Q}'_{Sp}${\it-action} if $f(k^*)\subseteq \tilde{Q}'_{Sp}$; these
are precisely the $Q'$-actions which preserve the space of symplectic
lattices.  Such a homomorphism induces a $Q'_{Sp}$-action via the
projection $GSp(V)\to PGSp(V)$ as described in section 2.

	We now assume that $\char\ k\ne 2$.  If $N\in\frak{sp}(V_F)$ is
semisimple, its eigenvalues must be of the form $\pm a_1,\dots,\pm a_n$
where $a_i\in\overline F$.  The characteristic polynomial
$\char_N(\mu)=\prod_{i=1}^n(\mu^2-a_i^2)$ is an element of $F[\mu^2]$, and
we define a polynomial $h_N\in F[\mu]$ by $h(\mu^2)=\char_N(\mu)$.  The
eigenvalues of $\text{ad } N:\frak g\to\frak g$ are $\pm a_i\pm a_j$ for
$i<j$, $\pm 2a_i$, and $0$ with multiplicity $n$.  The endomorphism $N$ is
regular semisimple if the kernel of $N$ has dimension $n$; this occurs when
the $a^2_i$'s are distinct and nonzero.  It is well-known that the
characteristic polynomial, or equivalently $h_N$, determines the
$Sp(V_{\overline F})$-conjugacy classes of regular semisimple elements of
$\frak{sp}(V_{\overline F})$.  Proposition 3.2 shows that the same is
true for regular semisimple classes in $\frak{sp}(V_F)$.

	Let $\eusm C$ be a regular semisimple conjugacy class in
$\frak{sp}(V_F)$, and suppose that there exists $N\in\eusm C$ and a
$Q'_{Sp}$-action $f$ almost commuting with $N$.  The endomorphism $N$ is
also regular semisimple viewed as an element of $\frak{gl}(V_F)$, so by
theorem 2.4, $\char_N$ is a homogeneous polynomial in $\mu$ and
$\pi^{q'}$ for some nonzero $q'\in\frac1{2n-1}\Z\cup\frac1{2n}\Z$.
Moreover, since $\det N\ne 0$, the $\mu^0\pi^{2nq'}$ term does not vanish,
implying that $2nq'\in\Z$.  It is immediate that $h_N$ is homogeneous in
$\mu$ and $\pi^q$ where $q=2q'\in\frac1n\Z$.  The converse can also be
proven by an argument similar to the proof of theorem 2.4.

\begin{theorem} Let $\eusm C$ be a regular semisimple conjugacy class
in $\frak{sp}(V_F)$.  Then there is a representative $N\in\eusm C$ almost
commuting with a $Q'_{Sp}$-action if and only if $h_{\eusm C}(\mu)$ is a
homogeneous polynomial in $\mu$ and $\pi^q$ for some nonzero $q\in\frac1n\Z$.
$\eusm C$ is nil if and only if $q>0$.
\end{theorem}
\noindent If we write $h_{\eusm
C}(\mu)=\sum_{i=0}^nc_i\pi^{qi}\mu^{n-i}$ and identify
$\frak{sp}(V_F)$ with $\sp(F)$ by choosing a symplectic basis for $V$, we can
define a representative $N=(a_{ij})$ for $\eusm C$ by:
\[a_{ij}=\begin{cases} 1 &\text{ if $i=j+1, j\in[1,n-1]$ or $(i,j)=(1,n+1)$}\\
-1 &\text{if $j=i+1, i\in[n+1,2n-1]$}\\
(-1)^ic_i\pi^{qi} &\text{if $i=j+n, j\in [1,n]$}\\
0 &\text{otherwise.}\end{cases} \]

\begin{corollary}  The regular semisimple conjugacy classes in
$\frak{sp}(V_F)$ obtained from almost commuting pairs $(N,f)$ with $f$ a
$Q'_{Sp}$-action are precisely those with $h_{\eusm
C}(\mu)=\prod_{i=1}^r(\mu^d-b_i\pi^s)$ where $rd=n$, $(s,d)=1$, $(\char'\
k,d)=1$, and the $b_i$'s are distinct elements of $k^*$.
\end{corollary}

\begin{proof}  If $\eusm C$ is such a class, then $h_{\eusm C}$ is a
homogeneous polynomial in $\mu$ and $\pi^q$.  Let $q=s/d$ with $(s,d)=1$ and
$d>0$.  All the eigenvalues are nonzero, so by corollary 2.3, each is
homogeneous of degree $q$.  If $a\pi^q$ is an eigenvalue, then it satisfies
the irreducible polynomial $\mu^d-a^d\pi^s$.  Consequently, every
irreducible factor of $h_{\eusm C}$ has degree $d$, implying that $n=rd$ for
some positive integer $r$.  Since these $r$ factors must be distinct, we
obtain the $r$ different $b_i$'s.  Finally, the factors are separable
polynomials if and only if $(\char'\ k,d)=1$.
\end{proof}

	We now find an explicit criterion for ellipticity of regular
semisimple conjugacy classes in $\frak{sp}(V_F)$.  Write $h_N=h_1\dots h_l$
where the $h_i$'s are distinct irreducible polynomials.  If we set
$h'_i(\mu)=h_i(\mu^2)$, then $\char_N=h'_1\dots h'_l$.  Let
$E_i=F[\mu]/(h_i)$; it is a field extension of degree $n_i=\deg h_i$ over
$F$.  By the Chinese remainder theorem,
$E=F[\mu]/(h_N)\cong\oplus_{i=1}^lE_i$.  Similarly, if we set
$E'_i=F[\mu]/(h'_i)$, $E'=F[\mu]/(\char_N)\cong\oplus_{i=1}^lE'_i$.  The
$F$-algebra injection $E\hookrightarrow E'$ given by $\mu\mapsto\mu^2$
induces injections $E_i\hookrightarrow E'_i$ for each $i$, making $E'_i$
into a $2$-dimensional $E_i$-algebra.  If $h'_i$ is irreducible, then
$E'_i$ is a field.  Otherwise, let $g_1,\dots,g_t$ be the monic,
irreducible factors of $h'_i$.  If $g_j(\mu)|h'_i$, then so does $g_j(-\mu)$.
Moreover, if the ideals generated by $g_j(\mu)$ and $g_j(-\mu)$ were equal, then
$g_j$ would be a polynomial in $\mu^2$ and would give rise to a proper
irreducible factor of $h_i$.  We can thus assume that
$g_2(\mu)=(-1)^ng_1(-\mu)$; the coefficient guarantees that $g_1$ and
$g_2$ are simultaneously monic.  But now $g_1g_2$ is a polynomial in
$\mu^2$, and this contradicts irreducibility of $h_i$ unless $h'_i=g_1g_2$.
Since $g_1$ and $g_2$ are distinct, $E'_i\cong F[\mu]/(g_1)\oplus
F[\mu]/(g_2)\cong E_i\oplus E_i$.  The last isomorphism follows because the homomorphisms
$E_i\to F[\mu]/(g_j)$ for $j\in[1,2]$ are nonzero maps between fields which have
the same dimension as $F$-vector spaces.

\begin{proposition}  Let $N$ be a regular semisimple element of
$\frak{sp}(V_F)$.  Then $N$ is elliptic if and only if $F[\mu]/(h_N)$ and
$F[\mu]/(\char_N)$ decompose into direct sums of the same number of fields.
Equivalently, for each irreducible factor $h_i$ of $h_N$,
$h'_i(\mu)=h_i(\mu^2)$ is also irreducible.
\end{proposition}

\begin{proof} We construct an element $N_i\in\frak{sp}_{2n_i}(F)$ with
$h_{N_i}=h_i$ by the procedure of theorem 4.2.  Using
the injection $\prod_{i=1}^l\frak{sp}_{2n_i}(F)\hookrightarrow\sp(F)$, we
obtain a `block-diagonal' matrix $N'\in\sp(F)$ with $h_{N'}=h_N$.  (The
restriction of $N'$ to $W_1=\oplus_{j=1}^{n_1}(Fe_j\oplus Fe_{n+j})$ is $N_1$
and so on.)  Since $N$ is regular semisimple, $N$ and $N'$ are conjugate,
and there is a symplectic basis $\underline w$ of $V_F$ in which $N$ is
given by $N'$.

	Suppose some $h'_i$ is reducible.  Since $V_F=W_1\perp\dots\perp W_l$, it
suffices to exhibit an $F$-cocharacter $F^*\to Z_{Sp(W_i)}(N_i)$.  We thus
assume without loss of generality that $h_N$ is irreducible, but $\char_N$
is not.  We then have $\char_N=g_1g_2$ with $g_2(\mu)=(-1)^ng_1(-\mu)$.
Construct a matrix $M\in\frak{gl}_n(F)$ with characteristic polynomial
$g_1$ as in the proof of theorem 2.4.  It has $n$ distinct
eigenvalues $a_1,\dots,a_n\in\overline F$ and $g_1(\mu)=\prod_{i=1}^n
(\mu-a_i)$.  The transpose $M^t$ has the same characteristic polynomial and eigenvalues,
so $\char_{-M^t}(\mu)=\prod_{i=1}^n (\mu+a_i)=(-1)^ng_1(-\mu)=g_2(\mu)$.
The block-diagonal matrix $M'=\text{diag}(M,-M^t)\in\sp(F)$ has the same
characteristic polynomial as $N$, and again invoking regularity, $N$ and
$M'$ are $\Sp(F)$-conjugate.  Since the $F$-cocharacter $F^*\to\Z_{\Sp(F)}(M')$
given by $x\mapsto\text{diag}(xI_n,x^{-1}I_n)$ is nontrivial,  $N$ is
not elliptic.

	Conversely, suppose $N$ is not elliptic, so that there exists a
nontrivial cocharacter $F^*\overset{\phi}{\to} Z_{\Sp(F)}(N)$.  Using the
above orthogonal decomposition of $V_F$, we obtain cocharacters
$F^*\overset{\phi_i}{\to} Z_{Sp(W_i)}(N_i)$.  Since at least one $\phi_i$
is nontrivial, we can assume without loss of generality that $h_N$ is
irreducible.  If $\char_N$ is also irreducible, then there cannot even be a
nontrivial cocharacter $F^*\to Z_{SL_{2n}(F)}(N)$.  Thus, $\char_N$ is
reducible as desired.
\end{proof}

	Suppose that $N$ is nil-elliptic and almost commutes with a
$Q'_{Sp}$-action.  The irreducible factors of $h_N$ are all of the form
$h_i(\mu)=\mu^d-b_i\pi^s$ with $b_i\in F^*$, $(d,s)=1$, and $s>0$.  The polynomial
$h'_i(\mu)=\mu^{2d}-b_i\pi^s$ is irreducible if and only if $(2d,s)=1$,
i.e. when $(d,s)=1$ and $s$ is odd.  Applying the previous result, we have:

\begin{proposition}  The nil-elliptic conjugacy classes in
$\frak{sp}(V_F)$ obtained from almost commuting pairs $(N,f)$ with $f$ a
$Q'_{Sp}$-action are precisely those with $h_{\eusm
C}(\mu)=\prod_{i=1}^r(\mu^d-b_i\pi^s)$ where $s>0$, $rd=n$, $(s,2d)=1$,
$(\char'\ k,d)=1$, and the $b_i$'s are distinct elements of $k^*$.  If $N$ is the
standard matrix representative of theorem 4.2, $f$ almost commutes with
exponent $s$ with the diagonal $Q'_{Sp}$-action $f$ given by
$\nu(m,i)=2dm+is$ and $\nu(m,n+i)=2dm+(1-i)s$ for $i\in [1,n]$.
\end{proposition}

	The varieties $(\hat{\eusm P}_N)^{k^*}$ are in general much simpler
than the $\hat{\eusm P}_N$'s, but, because the $k^*$-actions of the above
proposition do not always have discrete fixed point sets, they are not
necessarily finite. It is possible to identify the fixed lattices
precisely; the off-diagonal entries which may be nonzero as well as the
terms present in these Laurent series are fully determined by $r$.

	We now restrict ourselves to the case $r=1$, $d=n$, and $(\char'\ k,
2n)=(s,2n)=1$.  We let $N_{bs}$ denote the standard representative of the class with
$h_{\eusm C}(\mu)=\mu^n-b\pi^s$ given in theorem 4.2.  (This
does not seriously conflict with our previous notation in section 3
because the element of $\frak{sl}_{2n}(F)$ given by the two definitions are
in the same $SL_{2n}(F)$-conjugacy class.  Also, note that the fixed point
varieties for $N_{bs}$ viewed as an element of $\frak{sp}_{2n}(F)$ are closed
subvarieties of the analogous fixed point varieties with $N_{bs}$ viewed as
an element of $\frak{sl}_{2n}(F)$; in particular, they are algebraic
varieties even in positive characteristic.)  Again, we can assume $b=1$, and
we set $N_s\overset{\text{def}}{=}N_{1s}$. By the remark after lemma 3.4,
the $Q'_{Sp}$-action $f$ fixes only the diagonal lattices. The set
of symplectic lattice chains of type $J$ composed of diagonal lattices will
be denoted by $\Ds_J$; in particular, we set $\Ds_m=\eusm D_m\cap\Ls_m$.
Given $I\subseteq[1,n]$, let $\d^I_i=1$ if $i\in I$ and $\d^I_i=0$ otherwise.
It is immediate that for $0\le m\le n$, $\L_{\r}\in\Ds_m$ if and only if there
exists $I\subseteq[1,n]$ such that $|I|=m$ and $r_{n+i}= -r_i+\d^I_i$.

	We will again apply the theorem of Bialynicki-Birula to give
explicit combinatorial interpretations for the Euler
characteristics $\chi(\hat{\eusm P}_{N_{s}})$.  Because of the considerable
loss of symmetry involved in the passage to the symplectic algebras, the resulting counting problems are more complicated
than in the special linear case.

	We begin by determining the fixed point sets
$(\Ls_{m,N_{s}})^{k^*}$.  Suppose that $\L=\L_{\r}\in\Ds_m$ is fixed by $N_{s}$.  We must have $\pi^{r_i}e_{i+1}\in A\pi^{r_{i+1}}e_{i+1}$ and
$-\pi^{r_{n+i+1}}e_{n+i}\in A\pi^{r_{n+i}}e_{n+i}$ for $i\in[1,n-1]$,
$\pi^{r_{n+1}}e_1\in A\pi^{r_1}e_1$, and $(-1)^{n+1}\pi^{r_n+s}e_{2n}\in
A\pi^{r_{2n}}e_{2n}$.  As a result,
$(\Ls_{m,N_{s}})^{k^*}=\{\L_{\r}\mid\r\in\Rs\}$ where $\Rs$ is the
union over $I\subseteq[1,n]$ of size $m$ of the disjoint sets
\begin{multline*} \RsI=\{\r\in\Z^{2n}\mid r_n\le r_{n-1}\le\dots\le r_1\le r_{n+1}\le
r_{n+2}\le\dots\le r_{2n}\le r_n+s\\ \text{and $r_{n+i}=-r_i+\d^I_i$}\}.\end{multline*}

	Let $\R$ be the set of $2n$-tuples such that $(\eusm
L_{m,N_{s}})^{k^*}=\{\L_{\r}\mid\r\in\R\}$.  This is the same definition as
in (3-4) up to the numbering of the entries.  Thus, we have a bijection
$\R\overset{\phi}{\to}\C=\C^{2n}$ as in lemma 3.6.  The indexing
$\c=(c_1,\dots,c_{2n})$ is not convenient here, and we will write instead
\[(c_1,\dots,c_{2n})=(q_n,q_{n-1},\dots,q_1,q_0,q_{-1},\dots,q_{-(n-1)}).\]
The function $\phi$ is given by $q_0=r_{n+1}-r_1$, $q_n=r_n-r_{2n}+s$,
$q_i=r_i-r_{i+1}$, and $q_{-i}=r_{n+i+1}-r_{n+i}$ for $i\in[1,n-1]$.

	We now determine $\phi(\Rs)$.  Given a set $I\subseteq[1,n]$ and
$i\in[1,n-1]$, let $\e^I_i=\d^I_{i+1}-\d^I_i$.  Define \begin{multline*}
\CsI=\{\bold c\in C^{2n}_s\mid q_{-i}=q_i+\e^I_i\text{ for $i\in[1,n-1]$,
$q_n$ is even $\iff n\in I$,}\\ \text{and $q_0$ is odd $\iff 1\in
I$}\}\end{multline*} and \[\Cs=\bigcup_{\substack{I\subseteq[1,n]\\ |I|=m}}
\CsI.\] The union is clearly disjoint.  Note that it is not a priori
obvious that $\Cs\subseteq\C$.

	We will show that $\phi(\Rs)=\Cs$.  We need the following lemmas:

\begin{lemma}  For any $n\ge 1$ and any $I\subseteq[1,n]$ with
$|I|=m$, \[\sum_{i=1}^{n-1} i\e^I_i=\begin{cases} n-m& \text{if $n\in I$}\\ 
-m& \text{if $n\notin I$}.\end{cases}\] 
\end{lemma}

\begin{proof} We have \[\sum_{i=1}^{n-1} i\e^I_i=\sum_{i=1}^{n-1}
i(\d^I_{i+1}-\d^I_i)=n\d^I_n-\sum_{i=1}^n \d^I_i=n\d^I_n-m\] as desired.
\end{proof}

\begin{lemma}  $\Cs\subseteq\C$.
\end{lemma}

\begin{proof} We must show that $\sum_{i=1}^{2n}(l+i)c_i\equiv-m+(l+1)s\pmod
{2n}$ for some $l\in\Z$.  Set $l=-(n+1)$.  The sum is then
$-nq_n+0q_0+\sum_{i=1}^{n-1}i(q_{-i}-q_i)=-nq_n+\sum_{i=1}^{n-1}i\e^I_i$.
By lemma 4.6, this is congruent to $n-m$ regardless of the parity of $q_n$. 
Since $s$ is odd,  $-m-ns\equiv -m+n\pmod{2n}$, and the
result follows.  
\end{proof}

	We are now ready to prove:

\begin{proposition}  The map $\phi$ induces bijections
$\RsI\tilde{\to}\CsI$ for each $I\subseteq[1,n]$ and hence
bijections $\Rs\tilde{\to}\Cs$ for all $m\in[0,n]$.
\end{proposition}

\begin{proof} It is clear from lemma 4.7 that
$\phi(\RsI)\subseteq\CsI\subseteq C^{Sp}_{s|I|}$, so it will suffice to show
that $\phi^{-1}(\CsI)\subseteq \RsI$.  Choose $\c\in\CsI$, and let
$\r=\phi^{-1}(\c)$.  We will show by downward induction on $i$ that
$r_i+r_{n+i}=\d^I_i$.

	Using the formula for $\phi^{-1}$ given in lemma 3.6, we have
\begin{align*} r_{2n}&=\frac{m-s+\sum_{j=-n+1}^n (n+1-j)q_j}{2n}\quad\text{and}\\
r_n&=\frac{m-s+q_n+\sum_{j=-n+1}^{n-1} (-n+1-j)q_j}{2n}.\end{align*} 
Adding, we get
\begin{equation*}\begin{split} r_n+r_{2n}&=\frac{m-s+q_n+\sum_{j=-n+1}^{n-1}(1-j)q_j}{n}\\&=\frac{m-s+q_n+\sum_{j=-n+1}^{n-1}q_j-\sum_{j=-n+1}^{n-1}jq_j}{n}\\&=\frac{m-s+q_n+(s-q_n)+\sum_{j=1}^{n-1}
j\e^I_j}{n}=\d^I_n,\end{split}\end{equation*} 
where the last step is lemma 4.6.

	Now suppose $r_{i+1}+r_{n+i+1}=\d^I_{i+1}$.  By the formula for $\phi^{-1}$, $r_i+r_{n+i}=r_{i+1}+r_{n+i+1}+(q_i-q_{-i})$.  If $i+1\in I$, then
the term in parentheses is $-1$ (resp. $0$) for $i\notin I$ (resp. $i\in I$)
as desired.  The case $i+1\notin I$ is checked similarly.
\end{proof}

	Given a set $I\subseteq[1,n]$, define $\hat{I}=I\cup\{n+1\}$.  Let
\[\GsI=\{\p=(p_0,\dots,p_n)\in C_s^{n+1}\mid p_i\text{ is even $\iff i,i+1\in\hat{I}$ or $i,i+1\notin\hat{I}$}\},\] and let $\Gs$
denote the (disjoint) union of these sets over all $I$ of size $m$.  Define
a function $\CsI\overset{\psi}{\to}\GsI$ by $\q\mapsto
(q_0,q_1+q_{-1},\dots,q_{n-1}+q_{-(n-1)},q_n)$.  We also have a map
$\GsI\to\CsI$ with $q_0=p_0$, $q_n=p_n$, and for $i\in [1,n-1]$, 
\[(q_i,q_{-i})=\begin{cases} ([p_i/2],[p_i/2]+1)&\text{if $i\notin
I,i+1\in I$}\\([p_i/2]+1,[p_i/2])&\text{if $i\in I,i+1\notin
I$}\\([p_i/2],[p_i/2])&\text{otherwise.}\end{cases} \] 
It is easy to see that these maps are inverses, so we have proven:

\begin{lemma} The function $\psi$ defines bijections $\CsI\tilde{\to}\GsI$ and
$\Cs\tilde{\to}\Gs$ for all $I\subseteq[1,n]$ and $m\in[0,n]$.  
\end{lemma}

	For any $I\subseteq[1,n]$ with $|I|=m$, we can interpret $\GsI$ as
the set of distinct arrangement of $s$ balls in $n+1$ boxes with the parity
of the number of balls in each box determined by $I$.
Suppose $l$ of the boxes contain an odd number of balls.  Since these boxes
are nonempty and $s=2t+1$ is odd, $l=2a+1$ for some $a\in[0,t]$.  Note that
$\sum_{i=0}^n [p_i/2]=(s-l)/2=t-a$, and any such arrangement is uniquely
determined by the $[p_i/2]$'s.  Thus, there is a bijection from $\GsI$ to
$C^{n+1}_{t-a}$, implying that $|\GsI|=\binom{n+t-a}{n}$.

	Let $h_{n,m,a}$ be the number of subsets $I\subseteq[1,n]$ with
$|I|=m$ such that $2a+1$ of the $p_i$'s are odd.  A {\it succession} in
$\hat{I}$ is a pair of numbers $i,i+1\in\hat{I}$.  Let $g_{n,m,j}$ be the
number of $I$'s of size $m$ such that $\hat{I}$ has $j$ successions.  Write
$\hat{I}=\{i_1<\dots<i_m<i_{m+1}=n+1\}$, and partition it into maximal
sequences of consecutive integers $\{i_{l_{j-1}+1}<\dots<i_{l_j}\}$.
Suppose there are $d$ sequences in the partition. Since $0\notin\hat{I}$ and
$n+1\in\hat{I}$, $l_0=0$ and $l_d=m+1$.  Each maximal sequence contributes
$l_j-l_{j-1}-1$ successions to $\hat{I}$, so there are a total of
$l_d-l_0-d=m+1-d$ successions.  The boxes containing an odd number of balls are
those indexed by $i_{l_{j-1}+1}-1$ for $j\in[1,d]$ and $i_{l_j}$ for
$j\in[1,d-1]$.  Hence, there are $2(d-1)+1$ such boxes, i.e. $a=d-1$,
implying that $h_{n,m,a}=g_{n,m,m-a}$.  We must have $0\le a\le m$, since
the number of successions is evidently between $0$ and $m$.

\begin{lemma} For each $j\in [0,m]$, \[g_{n,m,j}=\binom{m}{j}\binom{n-m}{m-j}.\]
\end{lemma}

\begin{proof} Let $f_{n,m,j}$ be the number of size $m$ subsets
$I\subseteq[1,n]$ having $j$ successions.  It is well known (see
\cite[p. 12]{R}) that \[f_{n,m,j}=\binom{m-1}{j}\binom{n-m+1}{m-j}.\]  We
see that $g_{n,m,j}$ is
the number of subsets of $[1,n+1]$ of size $m+1$ which contain the element
$n+1$ and have $j$ successions.  This collection of subsets consists of those
size $m+1$ subsets with $j$ successions which are not already subsets of
$[1,n]$.  Thus, \begin{equation*}\begin{split}
g_{n,m,j}=f_{n+1,m+1,j}-f_{n,m+1,j}&=\binom{m}{j}\binom{n-m+1}{m+1-j}-\binom{m}{j}\binom{n-m}{m+1-j}\\
&=\binom{m}{j}\binom{n-m}{m-j}.\qed\end{split}\end{equation*} \renewcommand{\qed}{}\end{proof}

	We can finally compute $|\Gs|$.

\begin{proposition}
\[|\Gs|=\binom{[s/2]+m}{m}\binom{[s/2]+n-m}{n-m}.\]
\end{proposition}

\begin{proof} We have \begin{equation*}\begin{split} |\Gs|=\sum_{\substack{I\subseteq[1,n]\\ |I|=m}}
|\GsI|&=\sum_{a=0}^m h_{n,m,a}\binom{n+t-a}{n}=\sum_{a=0}^m
g_{n,m,m-a}\binom{n+t-a}{n}\\&=\sum_{a=0}^m
\binom{m}{m-a}\binom{n-m}{a}\binom{t+n-a}{n}\\&=\sum_{a=0}^m
\binom{m}{a}\binom{n-m}{a}\binom{t+n-a}{n}\\&=\binom{t+m}{m}\binom{t+n-m}{n-m}.\end{split}\end{equation*}
The last equality is proven in \cite[p. 17]{R}; it is an application of the
Vandermonde convolution formula.  The result follows since $t=[s/2]$.  Note
that the summation from $0$ to $m$ is correct even when $t$ is smaller than
$m$ because $\binom{n+t-a}{n}=0$ for $a>t$.  
\end{proof}

	Combining the fact that $\chi(\Ls_{m,N_{s}})=|\Rs|$ with lemma 4.9
and propositions 4.8 and 4.11, we obtain:

\begin{corollary}  For any $m\in[0,n]$,
\[\chi(\Ls_{m,N_{s}})=\binom{[s/2]+m}{m}\binom{[s/2]+n-m}{n-m}.\]
In particular, $\chi(\Ls_{m,N_{s}})=\chi(\Ls_{{n-m},N_{s}})$ and
$\chi(\Ls_{0,N_{s}})=\chi(\Ls_{n,N_{s}})=\binom{[s/2]+n}{n}$.
\end{corollary}

\begin{remarks} \begin{mrk}
\item The variety $(\Ls_{m,N_{s}})$ is isomorphic to
the varieties $(\Ls_{m+2ni,N_{s}})$ and $(\Ls_{-m+2ni,N_{s}})$ for any
$i\in\Z$ via the maps $L\mapsto\pi^i L$ and $L\mapsto\pi^i L^*$
respectively.  Thus, the result is true for all $m\in\Z$ if $m$ is replaced
by $m'\in[0,n]$ such that $m\equiv\pm m'\pmod{2n}$.
\item The formula for self-dual lattices can be proven in a much
simpler way.  Since $r_{n+i}=-r_i$ for all $i$, a self-dual lattice
$\L_{\bold r}$ is fixed by $N_{s}$ when $r_n\le r_{n-1}\le\dots\le r_1\le
-r_1\le\dots\le -r_n\le r_n+s$.  Setting $d_i=-r_i$, we have
$0=d_0\le d_1\le\dots\le d_n\le d_{n+1}=[s/2]$.  This is just an
arrangement of $[s/2]$ balls in $n+1$ boxes (the $i$th box contains
$d_i-d_{i-1}$ balls), and there are $\dbinom{[s/2]+n}{n}$ such
arrangements.
\end{mrk}
\end{remarks}  

\renewcommand{\Rs}{R^{Sp}}

	We now study the fixed symplectic lattice chains of type $J$ for
$J$ a nonempty subset of $[0,n]$.  Let $J=\{m_1<\dots< m_l\}$.  It
is clear from the definition that a lattice chain $(L_m)\in\Xs_J$ is
determined by the $L_{m_i}$'s.  Conversely, any sequence
$L_{m_1}\supset\dots\supset L_{m_l}$ with $L_{m_i}\in\Ls_{m_i}$ extends
uniquely to a chain in $\Xs_J$.  To see this, note that if $L\subset L'$,
then $L^*\supset L'{}^*$.  Thus, if $L$ and $L'$ are both symplectic, the
totally ordered set of lattices associated to $L$ and $L'$ interlock to
form a chain of lattices.  This shows that
$(\Xs_{J,N_{s}})^{k^*}=\{(\L_{\r^m})\mid
(\r^{m_1},\dots,\r^{m_l})\in\Ys_{Js}\}$ where
\[\Ys_{Js}=\{(\r^{m_1},\dots,\r^{m_l})\mid \r^{m_i}\in\Rs_{sm_i} \text{ and
$\r^{m_1}\le\dots\le\r^{m_l}$}\}.\] For $\r\in\Rs_{sm_1}$, let
$\Ys_{Js}(\r)$ denote the chains in $\Ys_{Js}$ with $\r^{m_1}=\r$.

	These sets can be given a uniform combinatorial description in
terms of a certain directed graph $\Delta$.  The vertex set of $\Delta$ is
$C^{n+1}_s$, 
the distinct arrangements of $s$ balls in $n+1$ boxes.  We index the boxes
from $0$ to $n$, and call the number of balls in the $i$th box $p_i$.
Arrange the boxes in a line with the index increasing to the right.
Viewing a box as two boundary walls, a vertex is just a string of $s$ balls
and $n$ walls, both numbered starting from $1$; the first and last walls
are redundant because they do not separate adjacent boxes.  Such a string
also uniquely determines an arrangement of $n$ walls in $s+1$ cells, where
a cell is determined by boundary balls.  Let $b_j$ be the number of walls
in the $j$th cell for $j\in[0,s]$, with the cell index also increasing to the
right.

	Two vertices $v$ and $v'$ are connected by an edge if one
configuration is obtained from the other by moving a ball between two
adjacent boxes, i.e. by reversing the order of a wall and ball next to each
other.  Given a vertex $v$, define $I_v\subseteq[1,n]$ recursively by
prescribing that for $i\in[0,n-1]$, the size of $\{i,i+1\}\cap I_v$ is congruent to
$p_i$ modulo $2$.
The recursion can begin because $0\notin I_v$.  If $v$ and $v'$ are adjacent,
then there exists $i\in[0,n-1]$ such that $p_j$ and $p'_j$ have the same
parity if and only if $j$ is neither $i$ nor $i+1$.  Thus, $I_v$ and $I_{v'}$ differ only
at the single element $i+1$.  We assign a direction to the edge from $v$ to
$v'$ by having it point towards the vertex with the larger associated set,
namely, the set containing $i+1$.

	The sets $\GsI$ are disjoint subsets of the vertex set and in fact
define a partition of $V(\Delta)$. 

\begin{proposition}  The vertex
set $V(\Delta)=C^{n+1}_s$ is the disjoint union of the $\GsI$'s.  Furthermore,
$\overrightarrow{vv'}$ is an edge if and only if there exists $m\in[0,n-1]$ such
that $\phi^{-1}\psi^{-1}(v)\in\Rs_{sm}$, $\phi^{-1}\psi^{-1}(v')\in
R^{Sp}_{s,m+1}$, and $\phi^{-1}\psi^{-1}(v)\le \phi^{-1}\psi^{-1}(v')$.
\end{proposition}

\begin{proof} Since the $\GsI$'s are disjoint, in order to prove the first claim, it
suffices to show that $v\in G^{Sp}_{sI_v}$.  It is
obvious from the construction of $I_v$ that $p_i$ has the parity required
of an element of $G^{Sp}_{sI_v}$ for $i\in[0,n-1]$.  It only remains to
check that $p_n$ is even if and only if $n\in I_v$.  Since $\sum_{i=0}^n
p_n=s$ is odd, this is a consequence of the following lemma.

\begin{lemma}  Suppose $\sum_{i=0}^j p_i$ is odd (resp. even).  Then $j\notin I_v$ if and only if $p_j$ is odd (resp. even).
\end{lemma}

\begin{proof} We use induction on $j$.  Since $0\notin I_v$, the statement
is vacuous for $j=0$.  Assume it is true for $j-1$.  Suppose $\sum_{i=0}^j
p_i$ is odd and $p_j$ is odd (resp. even).  Then $\sum_{i=0}^{j-1} p_i$ is even
(resp. odd).  If $p_{j-1}$ is even, then $j-1\notin I_v$ (resp. $j-1\in I_v$) by
induction.  Since $p_{j-1}$ even also implies that $j-1$ and $j$ are
simultaneously in or out of $I_v$, we have $j\notin I_v$ (resp. $j\in I_v$).  If
$p_{j-1}$ is odd, then $j-1\in I_v$ (resp. $j-1\notin I_v$), again giving
the desired result.  The proof is similar for $\sum_{i=0}^j p_i$ even.
\end{proof}

	Suppose that $\r\in\Rs_{sm}$, $\r'\in\Rs_{s,m+1}$, and $\r\le\r'$.
Then there exist subsets $I$ and $I'$ of $[1,n]$ such that $\r\in\Rs_{sI}$,
$\r'\in\Rs_{sI'}$, $|I|=m$, $|I'|=m+1$, and $I\subset I'$.  If $i$ is the
unique element of $I'\setminus I$, then $r'_j=r_j$
except at one of the two indices $i$ and $n+i$ where it is increased by
$1$.  Let $\p=\psi\phi(\r)\in G^{Sp}_{sI}$ and $\p'=\psi\phi(\r')\in
G^{Sp}_{sI'}$.  It is immediate that $\p'$ is obtained from $\p$ by shifting
a ball between the $(i-1)$th and $i$th boxes.  Since $I\subset
I'$, $\overrightarrow{\p\p'}$ is an edge.  The argument is reversible, so the
proposition is proved.  
\end{proof}

	Set $m_0=0$ and $m_{l+1}=n$.  Let $E^n_{Js}$ denote the set of
sequences of vertices $(v_1,\dots,v_l)$  such that $v_i\in G^{Sp}_{sm_i}$
and there exists a (directed) path in $\Delta$ going through all of the
$v_i$'s.  In particular, if $J=\{m,m+1,\dots,m'\}$, then $E^n_{Js}$ is the
set of paths from $\Gs$ to $G^{Sp}_{sm'}$.  For each $v\in G^{Sp}_{sm_i}$,
we let $E^n_{Js}(v)$ be the subset of $E^n_{Js}$ consisting of sequences
starting at $v$.  Proposition 4.13 implies that there is a
bijection $\Ys_{Js}\tilde{\to} E^n_{Js}$ which restricts to bijections $\Ys_{Js}(\r)\tilde{\to} E^n_{Js}(\psi\phi(\r))$.  We thus get a concrete
combinatorial interpretation of the Euler characteristics of the fixed
point varieties with respect to $N_{s}$ on each affine flag manifold:

\begin{proposition} Suppose that $\char\ k\ne 2$ and $(\char'\
k,n)=1$.  Let $s$ be an odd integer relatively prime to $n$ and $J$ a nonempty subset of $[0,n]$.  Then
$\chi(\Xs_{J,N_{s}})=|E^n_{Js}|$.  
\end{proposition}

	We define a set \[Z^n_{Js}=\{(v,\s)\mid \s:[1,n]\to [0,l],
\s^{-1}(0)=I_v, \text{and
$|\s^{-1}(i)|=j_i\overset{\text{def}}{=}m_{i+1}-m_i$}\}.\] Let
$\theta:E^n_{Js}\to Z^n_{Js}$ be the function
$(v_1,v_2,\dots,v_l)\mapsto(v_1,\s)$ where $\s^{-1}(0)=I_{v_1}$,
$\s^{-1}(l)=I_{v_l}^c$, and $\s^{-1}(i)=I_{v_{i+1}}-I_{v_i}$ for
$i\in[1,l-1]$.  (Note that $\s^{-1}(0)$ and $\s^{-1}(l)$ can be empty.)
Consequently, \begin{equation*}|E^n_{Js}|=\sum_{(v,\s)\in
Z^n_{Js}}\theta{}^{-1}((v,\s)).\tag{4-2}\end{equation*}

	This observation provides an algorithm for computing the Euler
characteristics in terms of a game.  If $v$ is a vertex viewed as a string
of balls and walls, then we obtain all $v'$ adjacent to $v$ by replacing
some substring $|\bullet$ by $\bullet|$ or vice versa, where $\bullet$
denotes a ball and $|$ a wall.  The edge $\overrightarrow{vv'}$ is contained in $\Delta$ if and only if
the wall involved in the switch is not in $I_v$.  The elements of
$\theta{}^{-1}((v,\s))$ may be found by the following procedure.  Find all
possible configurations of balls and walls obtained from $v$ by one move for each
wall in $\s^{-1}(1)$.  For each such $v_2$, repeat the process with the
walls in $\s^{-1}(2)$.  We continue until a ball has been moved over each
wall in $\s^{-1}([1,l-1])$, thereby getting all $l$-tuples of vertices
$(v,v_2,\dots,v_l)\in\theta{}^{-1}((v,\s))$.

\renewcommand{\a}{\bold a}

	In order to evaluate the sum (4-2), we will pass from $Z^n_{Js}$
to a more familiar index set consisting of intersection matrices.  As a
first step, we define a map $\eta:V(\Delta)\to C_n^{t+1}$ where $t=[s/2]$.
Given a vertex $v$, let $\b\in C^{s+1}_n$ be the corresponding arrangement of $n$ walls
in $s+1$ cells, and set $\eta(v)=\a$ with $a_k=b_{2k}+b_{2k+1}$.  This
function simply pools the walls in cells $2k$ and $2k+1$ into one larger
cell.  We let $\eta_m$ denote the restriction of $\eta$ to $\Gs$.  We also
set $C_{m}^{t+1}(\a)=\{(c_0,\dots,c_t)\mid c_k\in[0,a_k]\text{ and
$\sum_{k=0}^t c_k=m$}\}$.

\begin{proposition} There is a bijection between
$\eta_m^{-1}(\a)$ and $C_{m}^{t+1}(\a)$; in particular, $\eta_0$ is a
bijection.  Moreover, identifying $\Go$ and $C_n^{t+1}$ via $\eta_0$,
$\eta(v)$ is the unique element of $\Go$ from which $v$ can be obtained
using only moves of the form $|\bullet\rightsquigarrow\bullet|$.
\end{proposition}

\begin{proof} Suppose that $\p\in \Go$, so that each box contains an even
number of balls except for the last one.  The $i$th wall has
$p_0+\dots+p_{i-1}$ balls to its left, and  since this is an even number,
it is in a cell with an even index.  Thus, $b_j=0$ for $j$ odd.  In this
case, $a_k=b_{2k}$ and $\eta_0$ is injective because $\b$ is uniquely
determined by $\p$.  Surjectivity follows since $|\Go|=|C_n^{t+1}|$ by
proposition 4.11.

	Consider a game starting at $\p$ with only moves of the form
$|\bullet\rightsquigarrow\bullet|$ allowed.  The only balls which can take
part in these switches are the leftmost balls in boxes $1,\dots,n$.  To see
this, note that walls can only move to the right.  Consequently, the balls
in box $0$ play no role.  Also, if a ball becomes the first ball in its box
during the course of the game without starting out in this position, then
the wall left of it must already have moved.  Since each wall can only move
once, the ball is not involved in the game.  It follows that no even-numbered balls can move.  Thus, $b_{2k}+b_{2k+1}$, the number of walls
between balls $2k$ and $2(k+1)$, never changes, i.e. $\eta$ is constant on
the set of configurations obtained as outcomes of this game.  This proves the uniqueness part of the final
statement.

	We will show existence of such a path by induction on $m$.  This is
trivial for $m=0$.  Suppose $m\ge1$ and the claim is true for $m-1$.  Take
$v=\p\in\Gs$, and let $i$ be the smallest integer such that $p_i$ is odd.
Equivalently, $i+1$ is the smallest element of $\hat{I}_v$; in fact,
$i+1\in I_v$ because $m\ne
0$.  Since $p_i\ne 0$, the $(i+1)$st wall is immediately
preceded by a ball.  We obtain a vertex $v'$ adjacent to $v$ by moving this
wall to the left, and $\overrightarrow{v'v}$ is an edge since $i+1\in I_v$.
This implies that, $v'\in G^{Sp}_{s,m-1}$, so there exists a path to $v'$ of the correct
form by inductive hypothesis.  Concatenating $\overrightarrow{v'v}$ with this
path proves the result.

	Finally, given $\c\in C_{m}^{t+1}(\a)$, we obtain a vertex
$\tau(\c)$ in $\Gs$ by moving ball $2k+1$ past the last $c_k$ walls in the
$2k$th cell; such a
ball exists because the last cell has odd index.  These moves are of the
specified type, so $\eta(\tau(\c))=\a$.   The map
$C_{nm}^{t+1}(\a)\overset{\tau}{\to}\eta_m^{-1}(\a)$  is
surjective by the previous paragraph.  It is injective because the $(2k+1)$st
cell of $\tau(\c)$ contains $c_k$ walls.
\end{proof}

	Let $\bold m_J$ be the sequence $(0=m_0,m_1,\dots,m_l,m_{l+1}=n)$.  We define a map
$\xi$ from $Z^n_{Js}$ to the set of $(\m_J,t+1)$-intersection matrices
$Q_{\m_J,t+1}$.  (Both row and column indexing start at $0$.)  Let
$\a=\eta(v)$, and let $A_k$ be the set of walls in the $2k$th cell of
$\eta_0^{-1}(\a)$.  Set $\xi(v,\s)=(|\s^{-1}(i)\cap A_k|)$ and
$\zeta=\xi\circ\theta$.  We obtain \begin{equation*} |E^n_{Js}|=\sum_{Q\in
Q_{\m_J,t+1}}|\zeta^{-1}(Q)|.\tag{4-3}\end{equation*}

	Fix an intersection matrix $Q$.	Notice that the only balls which play a
role in moves involving the walls in $A_k$ are the boundaries of the $2k$th
cell.  Since a ball bounds only one even-numbered cell, this shows that the
manipulations of the game take place independently within each
$A_k$.  Consequently, it will suffice to understand the game for the
special cases of a string of walls with a boundary ball on both sides or
only on the right.  These correspond to the game restricted to the $2k$th
cell for $k>0$ and $k=0$ respectively.

	Assume the cell contains $d$ walls and suppose $z_0,\dots,z_l$ are
nonnegative integers such that $\sum_{i=1}^lz_i=d$.  Consider a game with $z_i$ moves at the $i$th
step for $i\in[0,l-1]$, where only moves of the form
$|\bullet\rightsquigarrow\bullet|$ are allowed at the $0$th step.  Call the
number of outcomes of the game $\alpha_{\bold z}^d$ if the cell is only
bounded on the right and $\beta_{\bold z}^d$
if it is bounded on both sides.

\begin{lemma} \[\alpha_{\bold z}^d=1 \text{ and $\beta^d_{\bold
z}=\prod_{i=1}^{l-1} (z_i+1)$}.\] 
\end{lemma}

\begin{proof} Suppose there is only a boundary ball on the right.  Then
there is only one outcome of the game; indeed, we have already seen in the
proof of proposition 4.16 that each move consists of the ball shifting one
box to the left.  Now suppose the cell has a ball on each side.  The $0$th
move is uniquely determined, and setting $z'_0=0$ and $z'_i=z_i$ for
$i\in[1,l]$, we have $\beta^d_{\bold z}=\beta^{d-z_0}_{\bold z'}$.  We may thus assume
that $z_0=0$.  The left ball can only move to the right and vice versa
because otherwise a ball would jump over the same wall twice.  Hence, for
each $c\in[0,d]$, there are $c+1$ possible configurations after $c$ moves
corresponding to $c'\in[0,c]$ jumps of the left ball and $c-c'$ jumps of
the right ball.  The resulting strings have the form of the following
diagram:
\[\overbrace{||\cdots||}^{c'}\bullet\overbrace{||\cdots||}^{d-c}\bullet\overbrace{||\cdots||}^{c-c'}.\]
The outer walls have all been used, so the $c$ moves reduce the game to the
game of the same form with $d-c$ walls.  The lemma follows from these
observations by induction.  
\end{proof}

\begin{corollary}
\[|\zeta^{-1}(Q)|=\prod_{i=1}^{l-1}\prod_{k=1}^t (q_{ik}+1)\]
\end{corollary}

\begin{proof}  Let $\bold q^k$ be the $k$th column of $Q$.  Then
$|\zeta^{-1}(Q)|=\alpha^{a_0}_{\bold
q^0}\prod_{k=1}^t\beta^{a_k}_{\bold q^k}$, and lemma 4.17 completes the
proof.
\end{proof}

	  Let $\gamma_d^t=\sum_{\bold y\in C^t_d}\prod_{k=1}^t (y_k+1)$.
We need the following lemma:

\begin{lemma} For any $d\ge 0$ and $t>0$,
\[\gamma_d^t=\binom{2t+d-1}d.\]
\end{lemma}
\begin{proof} We use induction on $t$.  The claim is true for $t=1$, since
$\gamma_d^1=d+1=\binom{d+1}d$.  Now take $t>1$ and assume that the result holds
for $t-1$.  Then \begin{equation*}\begin{split} \gamma_d^t=\sum_{\bold y\in C^t_d}\prod_{k=1}^t
(y_k+1)&=\sum_{y_d=0}^d(y_d+1)\sum_{\bold x\in
C^{t-1}_{d-y_d}}\prod_{k=1}^{t-1}
(x_k+1)\\&=\sum_{y_d=0}^d(y_d+1)\binom{2(t-1)+d-y_d-1}{d-y_d}\\&=\sum_{y_d=0}^d\binom{y_d+1}{y_d}\binom{2(t-1)+d-y_d-1}{d-y_d}\\&=\binom{2t+d-1}d,
\end{split}\end{equation*}
where the last step is a version of the Vandermonde convolution formula \cite[p. 10]{R}.
\end{proof}

\begin{corollary} For any $d\ge 0$ and $t>0$, \[\sum_{\bold z\in C_d^{t+1}}\prod_{k=1}^t (z_k+1)=\binom{2t+d}d.\]
\end{corollary}
\begin{proof}  We have \[\sum_{\bold z\in C_d^{t+1}}\prod_{k=1}^t
(z_k+1)=\sum_{z_d=0}^d\gamma_{d-z_d}^t=\sum_{z_d=0}^d\binom{2t+d-z_d-1}{d-z_d}=\binom{2t+d}d;\]
the last equality is a standard binomial identity \cite[p. 7]{R}.
\end{proof}

	We can now calculate $|E^n_{Js}|$.

\begin{proposition}
\[|E^n_{Js}|=\binom{t+j_0}{t}\binom{t+j_l}{t}\prod_{i=1}^{l-1}\binom{s+j_i-1}{j_i}.\]
\end{proposition}
\begin{proof} Let $\bold q_i$ be the $i$th row of $Q$.  Recall
from lemma 3.9 that the map $Q\mapsto(\bold q_0,\dots,\bold q_l)$
is a bijection of $Q_{\bold m_J,s}$ onto $\prod_{i=0}^l C^{t+1}_{j_i}$.
From (4-3) and corollaries 4.18 and 4.20, we have \begin{equation*}\begin{split} |E^n_{Js}|&=\sum_{Q\in
Q_{\m_J,t+1}}|\zeta^{-1}(Q)|
=\sum_{\substack{(\bold q_0,\dots,\bold
q_l)\in\\ \prod_{i=0}^l C^{t+1}_{j_i}}}\ \prod_{i=1}^{l-1}\prod_{k=1}^t
(q_{ik}+1)\\
&=|C^{t+1}_{j_0}||C^{t+1}_{j_l}|\prod_{i=1}^{l-1}\sum_{\bold q_i\in
C^{t+1}_{j_i}}\prod_{k=1}^t(q_{ik}+1)=\binom{t+j_0}{t}\binom{t+j_l}{t}\prod_{i=1}^{l-1}\g_{j_i}^{t+1}\\&=\binom{t+j_0}{t}\binom{t+j_l}{t}\prod_{i=1}^{l-1}\binom{2t+j_i}{j_i}=\binom{t+j_0}{t}\binom{t+j_l}{t}\prod_{i=1}^{l-1}\binom{s+j_i-1}{j_i}.
\qed\end{split}\end{equation*}\renewcommand{\qed}{}
\end{proof}
\noindent Combining this result with proposition 4.15 gives the following
theorem:

\begin{theorem} Suppose that $\char\ k\ne 2$ and $(\char'\ k,n)=1$.
Let $s$ be an odd integer relatively prime to $n$ with $t=[s/2]$ and $J=\{m_1<\dots <m_l\}$ a nonempty subset of $[0,n]$.
Let $j_i=m_{i+1}-m_i$ where $m_0=0$ and $m_{l+1}=n$.  Then
\[\chi(\Xs_{J,N_{s}})=\binom{t+j_0}{t}\binom{t+j_l}{t}\prod_{i=1}^{l-1}\binom{s+j_i-1}{j_i}.\]
In particular, $\chi(\hat{\eusm B}_{N_{s}})=s^n$ and
$\chi(\Ls_{m,N_{s}})=\dbinom{t+j_0}{t}\dbinom{t+j_l}{t}$.  
\end{theorem}

\begin{proof}  It only remains to check the final statement.  The symplectic
lattices of type $m$ correspond to $J=\{m\}$, so $l=1$, $j_0=m$, and
$j_1=n-m$.  Substituting, we obtain another proof of corollary 4.12.  The
full affine flag manifold is isomorphic to $\Xs_{[0,n]}$.  In this case, $l=n+1$, $j_0=j_{n+1}=0$,
and $j_i=1$ for $i\in[1,n]$, and the result follows.
\end{proof}

\begin{remark}  As in the case of the root system $A_n$, the Euler characteristic is just the number of
simplices of the type corresponding to the parahoric class within the
subset of the positive Weyl chamber given by $\lambda_0\le s$, where
$\lambda_0$ is the affine root.
\end{remark}

\begin{corollary} Under the hypotheses of the previous theorem,
\[\chi(\Xs_{J,N_{s}})\le \chi(X_{J\cup(2n-J)\setminus\{2n\},N_{s}}),\]
with equality for $n=1$ or $s=1$ and strict inequality otherwise.
\end{corollary}

\begin{proof} The inequality follows immediately from the fact that
$(\Xs_{J,N_{s}})^{k^*}$ is a subset of the finite set
$(X_{J\cup(2n-J)\setminus\{2n\},N_{s}})^{k^*}$.  These sets are equal when
$n=1$ because $Sp_2(k)=SL_2(k)$; thus the formulas for the
Euler characteristics of the affine flag varieties of $Sp_2(F)$ given in
theorem 4.22 agree with those given in theorem 3.13.  Also, if $s=1$,
then the Euler characteristics are $1$ for both the symplectic and special
linear groups.

Now suppose that $n$ and $s$ are both greater than $1$.  We have
\[\chi(X_{J\cup(2n-J)\setminus\{2n\},N_{s}})=s^{-1}\binom{2t+2j_0}{2t}\binom{2t+2j_l}{2t}\prod_{i=1}^{l-1}\binom{s+j_i-1}{j_i}^2.\]
 Note that $\binom{2t+2q}{2t}=\binom{t+q}t$ for $q=0$ and 
\[\binom{2t+2q}{2t}=\binom{t+q}t\prod_{p=1}^q\frac{2t+2p-1}{2p-1}\ge
s\binom{t+q}t\] for $q\ge 1$.  Since $t\ge 1$, the inequality is
strict when $q\ge 2$.  Furthermore, $\binom{s+q-1}q^2\ge
s\binom{s+q-1}q$ for $q\ge 1$, with strict equality for $q\ge 2$.  This
proves the claim if $j_i\ge 2$ for any $i\in[0,l]$.  Otherwise, $\sum_{i=0}^l
j_i=n\ge 2$ implies that there exist two indices $i_1$ and $i_2$ such that
$j_{i_1}=j_{i_2}=1$.  The result follows because $s^2>s$.
\end{proof}

	We can also obtain a combinatorial description of the Euler
characteristics of the fixed point varieties in the classical symplectic
flag manifolds for certain nilpotent conjugacy classes of $\sp(k)$.  The
nilpotent conjugacy classes of $\sp(k)$ are characterized by symplectic
partitions of $2n$, i.e. partitions whose odd parts come in pairs
\cite{Sp}.  These partitions are the dimensions of the Jordan blocks of the
nilpotent element of $\sp(k)$.   Let $\eusm C$ be a nilpotent conjugacy
class corresponding to a partition with at most one part appearing with odd
multiplicity.  Such a class is determined by $n_0\in[0,n]$ and a partition
$n_1,\dots,n_l$ of $n-n_0$.  Let $M_0$ and $M_i$ be matrices with a single
Jordan block in $\frak{sp}_{2n_0}(k)$ and $\frak{gl}_{n_i}(k)$
respectively.  Recall that $\frak{gl}_n(k)$ injects into $\sp(k)$ via
$M\mapsto\text{diag}(M,-M^t)$.  Then the image of $(M_0,M_1,\dots,M_l)$
under the injections
$\frak{sp}_{2n_0}(k)\times\prod_{i=1}^t\frak{gl}_{n_i}(k)\hookrightarrow\prod_{i=0}^l\frak{sp}_{2n_i}(k)\hookrightarrow
\sp(k)$ is an element in this class; it has two Jordan blocks of size $n_i$
and one of size $2n_0$.

\begin{theorem} Assume that $\char\ k\ne 2$ and $(n,\char'\ k)=1$.
Let $N\in\sp(k)$ be an element of a nilpotent conjugacy class with at
most one Jordan block size having odd multiplicity, and let
$2n=2n_0+2\sum_{i=1}^l n_i$ be the corresponding symplectic partition of $2n$.
Let $s=2t+1$ be any odd integer greater than $2l$ and relatively prime to
$n$.  Define $\a\in C^{t+1}_n$ by $a_k=n_k$ for $k\in[0,l]$ and $a_k=0$ for $k>l$.
Then for each nonempty $J\subseteq[1,n]$, \[\chi(\eusm
F^{Sp}_{J,N})=|E^n_{J\cup\{0\},s}(\eta^{-1}(\a))|=\sum_{Q\in
Q_{\m_J,t+1}(\a)}\prod_{i=1}^{l-1}\prod_{k=1}^t (q_{ik}+1).\] In
particular, \[\chi(\eusm B^{Sp}_{N})=\frac{n!\,2^{n_1+\dots+
n_l}}{n_0!\cdots n_l!}.\] \end{theorem}

\begin{proof} Let $\r\in R^{Sp}_{s0}$ correspond to $\a$, so that setting
$\q=\phi(\r)$ and $\p=\psi(\q)$, we have $\eta_0(\p)=\a$.  We show that the
nilpotent map $\bar{N}_{s}$ induced on $W=\L_{\r}/\pi\L_{\r}$ has Jordan block
structure given by $n_0,\dots,n_l$.  It is clear from the definition that
$\eta_0^{-1}(\a)$ is the string of $s+1$ balls and $n$ walls with $n_k$
walls in the $2k$th cell for $k\in[0,l]$ and the other cells empty.  Thus,
the nonempty boxes are precisely those numbered $\sum_{i=0}^d n_i$ for
$d\in[0,l]$.  The nonzero $q_i$'s then have indices $n$ and
$\pm\sum_{i=0}^d n_i$ for $d\in[0,l-1]$.  As shown in the proof of theorem
3.11, a Jordan block of size $d$ corresponds to $j\in[-(n-1),n]$ such
that $q_j=q_{j+d}=0$ and $q_{j+d'}\ge 1$ for $0<d'<d$.  (It is not
necessary to understand the addition of the indices modulo $2n$ because
$q_n\ne 0$.)  It is immediate that there is one block of size $2n_0$ and
two of size $n_i$ for each $i\in[1,l]$ and that the block structure comes
from an injection of the form
$\frak{sp}_{2n_0}(k)\times\prod_{i=1}^t\frak{gl}_{n_i}(k)\hookrightarrow
\sp(k)$ as described above.  The first statement follows, since \begin{equation*}\begin{split}
\chi(\eusm F^{Sp}_{J,N})=|E^n_{J\cup\{0\},s}(\eta^{-1}(\a))|&=\sum_{Q\in
Q_{\m_J,t+1}(\a)}|\zeta^{-1}(Q)|\\ &=\sum_{Q\in
Q_{\m_J,t+1}(\a)}\prod_{i=1}^{l-1}\prod_{k=1}^t (q_{ik}+1). \end{split}\end{equation*}

	Now suppose $J=[1,n]$.  A $(\bold m_{[0,n]},s)$-intersection matrix
has a single $1$ in each row, so there are $\frac{n!}{a_0!\cdots a_t!}$
elements $Q\in Q_{\bold m_{[0,n]},s}(\a)$.  Also,
$|\zeta^{-1}(Q)|=2^{a_1+\dots+a_t}$, and since this is constant on the fiber over $\a$, $\chi(\eusm B^{Sp}_{N})=\dfrac{n!\,2^{a_1+\dots
+a_t}}{a_0!\cdots a_t!}=\dfrac{n!\,2^{n_1+\dots+
n_l}}{n_0!\cdots n_l!}$.  
\end{proof}

\begin{remarks} \begin{mrk} \item Since $\dfrac{n!\,2^{a_1+\dots a_t}}{a_0!\cdots a_t!}$ is the term corresponding to $\a$ in the multinomial expansion of
$s^n=(1+2+\dots+2)^n$, we obtain another proof of the formula for the
Euler characteristic of the full affine flag manifold.
\item  The Euler characteristic of each $\eusm P_N$ can be
obtained from the game with $2n+1$ balls.
\end{mrk}
\end{remarks}

\begin{corollary}  Assume $\char\ k\ne 2$ and $(n,\char'\ k)=1$.  Let
$N\in\sp(k)$ be an element of a nilpotent conjugacy class with at
most one Jordan block size having odd multiplicity.  Then for each nonempty $J\subseteq[1,n]$,
$\chi(\eusm F^{Sp}_{J,N})\le\chi(\eusm F_{J\cup(2n-J)})$.
\end{corollary}

\begin{proof}  The proof is similar to the first part of corollary 4.23.
\end{proof}

\begin{remark}  The formulas given in theorems 4.24 and
3.11 coincide for $\dim V=2$.
\end{remark}

\end{document}